\documentclass{article}

\usepackage{amsfonts}
\usepackage{amsmath,amsthm}

\usepackage[all]{xy}
\usepackage{graphicx}

\usepackage{amssymb}
\usepackage{latexsym}

\DeclareMathAlphabet\EuFrak{U}{euf}{m}{n}	
\SetMathAlphabet\EuFrak{bold}{U}{euf}{b}{n}	

\newcommand{\ra}{\rightarrow}
\newcommand{\lra}{\longrightarrow}

\newcommand{\hra}{\hookrightarrow}

\newcommand{\ovl}{\overline}

\newcommand{\wa}{\widehat}
\newcommand{\wt}{\widetilde}

\newcommand{\sC}{{\it C*}-}
\newcommand{\bC}{{\mathbb C}}
\newcommand{\bR}{{\mathbb R}}
\newcommand{\bT}{{\mathbb T}}

\newcommand{\bZ}{{\mathbb Z}}
\newcommand{\bM}{{\mathbb M}}
\newcommand{\bN}{{\mathbb N}}
\newcommand{\bP}{{\mathbb P}}

\newcommand{\bO}{{\mathbb O}}

\newcommand{\sud}{{{\mathbb {SU}}(d)}}

\newcommand{\eps}{\varepsilon}

\newcommand{\mA}{\mathcal A}
\newcommand{\mB}{\mathcal B}
\newcommand{\mC}{\mathcal C}
\newcommand{\mE}{\mathcal E}
\newcommand{\mF}{\mathcal F}
\newcommand{\mG}{\mathcal G}
\newcommand{\mH}{\mathcal H}

\newcommand{\mL}{\mathcal L}
\newcommand{\mM}{\mathcal M}

\newcommand{\mO}{\mathcal O}

\newcommand{\mZ}{\mathcal Z}

\newcommand{\aut}{ {\bf aut} }

\newcommand{\mSG}{ S_X(\mG) } 
\newcommand{\mSY}{ S_X(Y) }

\newcommand{\HG}{\mH \backslash \mG}
\newcommand{\HsG}{\mH_s \backslash \mG}

\newcommand{\goa}{ C(\mG) \otimes_X \mA }
\newcommand{\gof}{ \mF_s^\mG }
\newcommand{\gbof}{ S_\beta ( \mG , \wa \mF_s ) }

\newtheorem{thm}{Theorem}[section]
\newtheorem{cor}[thm]{Corollary}
\newtheorem{lem}[thm]{Lemma}
\newtheorem{prop}[thm]{Proposition}

\newtheorem{defn}[thm]{Definition}

\theoremstyle{definition}
\newtheorem{ex}{Example}[section]

\theoremstyle{remark}
\newtheorem{rem}{Remark}[section]

\numberwithin{equation}{section}

\begin{document}

\author{{\sf Ezio Vasselli\footnote{The author was partially supported by the European Network ``Quantum Spaces - Noncommutative Geometry" HPRN-CT-2002-00280.}
}
                         \\{\it Dipartimento di Matematica}
                         \\{\it University of Rome "La Sapienza"}
			 \\{\it P.le Aldo Moro, 2 - 00185 Roma - Italy }
                         \\{\sf vasselli@mat.uniroma2.it}}

\title{ Some Remarks on Group Bundles and \sC dynamical systems }
\maketitle

\begin{abstract}
We introduce the notion of fibred action of a group bundle on a $C(X)$-algebra. By using such a notion, a characterization in terms of induced \sC bundles is given for \sC dynamical systems such that the relative commutant of the fixed-point \sC algebra is minimal (i.e., it is generated by the centre of the given \sC algebra and the centre of the fixed-point \sC algebra). A class of examples in the setting of the Cuntz algebra is given, and connections with superselection structures with nontrivial centre are discussed.

\

\noindent {\em AMS Subj. Class.: 46L05, 55R10, 22D45}

\end{abstract}


\section{Introduction.}
\label{intro}

A result by S. Doplicher and J.E. Roberts (\cite[Thm.1]{DR86}) characterizes any compact \sC dynamical system $(\mB , G)$ with fixed-point algebra $\mA := \mB^G$ and conditions $\mA' \cap \mB = \mB \cap \mB'$, $\mA \cap \mA' = \bC 1$, in terms of an induced bundle of \sC algebras over a homogeneous space. The fibre of such a bundle is a \sC dynamical system $(\mF , H)$ with $H \subseteq G$, satisfying the conditions $\mF^H \simeq \mA$, $\mA' \cap \mF = \bC 1$.

It is of interest to generalize the above theorem to the case in which the centre of $\mA$ is non-trivial. From a mathematical point of view, our motivation is a duality theory for (noncompact) groups acting on Hilbert bimodules (\cite{Vas05}). 
From the point of view of mathematical physics, the original motivation was the following: given our \sC dynamical system $( \mB , G )$, define $\mC :=$ $\mB' \cap \mB$, $\mZ :=$ $\mA' \cap \mA$, and assume that $\mA' \cap \mB$ is generated as a \sC algebra by $\mZ$ and $\mC$; we look for a \sC dynamical system $( \mF , H )$, $H \subseteq G$, such that
\begin{equation}
\label{eq_rc}
\mF^H = \mA  \quad , \quad  \mA' \cap \mF = \mZ \ \ .
\end{equation}
\noindent The above relations play an important role in the context of certain \sC dynamical systems arising from superselection structures with nontrivial centre, called "Hilbert \sC systems" (see \cite{BL03}). From a physical viewpoint, (\ref{eq_rc}) appeared in low-dimensional conformal quantum field theory as a principle (see \cite[\S 1, p.142]{MS90}): the \sC algebra $\mF$ plays the role of a field algebra, $\mA$ is a universal algebra from which the "real" observable algebra can be recovered, and $H$ is the "gauge group", which in the above reference turns out to be actually a {\em quantum} group. In particular, the above interpretation of (\ref{eq_rc}) has been recently recognized (\cite[(4.34)]{Cio03}) in the Streater-Wilde model (\cite{SW71}).

Now, at the algebraic level $\mF$ appears as the image of a \sC epimorphism $\eta : \mB \ra \mF$, and $H$ is the stabilizer of $\ker \eta$ w.r.t. the $G$-action. In the case in which the centre of $\mA$ is nontrivial, it is easy to construct examples such that $H$ reduces to the trivial group (see Sec.\ref{ss_cuntz}). 
This unpleasant fact led us to introduce the following construction. As first, it is well-known that if an Abelian \sC algebra $C(X)$ is a unital subalgebra of the centre of a \sC algebra $\mB$, then we may regard $\mB$ as a bundle of \sC algebras with base space $X$. Given a bundle $\mG \ra X$ of compact groups, we introduce the notion of {\em fibred action} of $\mG$ on $\mB$, in such a way that each fibre of $\mG$ acts on the corresponding fibre of $\mB$ (Def.\ref{fibered_act}). 

Then, we adopt this point of view for the case of our \sC dynamical system $( \mB , G )$, by choosing $C(X) :=$ $\mC \cap \mA$. In this way, we find that $\mF$ has the structure of a bundle with base space $X$, on which a bundle $\mH$ acts in the above sense, with $\mH \subseteq$ $X \times G$. The fibred action of $\mH$ turns out to be non-trivial also in cases in which the stabilizer $H$ is trivial (Rem.\ref{rem_fa_q}).
Roughly speaking, we obtain that $\mF$ is a "field algebra" carrying an action by a bundle $\mH$ of "local gauge groups" (Thm.\ref{fix_point_3}, Thm.\ref{cor_ds}). 

In analogy with the original motivation of Doplicher and Roberts, we apply our construction to symmetric endomorphisms arising from superselection structures with non-trivial centre (Thm.\ref{thm_dr}, Thm.\ref{thm_bl}); in this setting, $C(X)$ is interpreted as the \sC subalgebra of $\mA' \cap \mA$ which is left invariant by the action of such endomorphisms.
Since in low-dimensional quantum field theory unitary braidings arise, we hope to develop in a future work an analogous construction, involving Hopf $C(X)$-algebras (\cite{Bla96}) instead of group bundles.

Groups of sections (of trivial bundles) acting on \sC algebras have been also considered in \cite{CG04}. In this case, the main motivation arises from local quantum gauge field theory, and $X$ is interpreted as a (compactified) space-time.

The present work is organized as follows. 

In Sec.\ref{bundles}, we give some basic properties about bundles $p : Y \ra X$, regarded as the topological counterparts of the commutative $C_0(X)$-algebras $C_0(Y)$. 
Then, we consider group bundles, and the associated groups of sections (for brevity called {\em section groups}, see Def.\ref{def_sg}); in the case in which the group bundle is locally trivial, we prove the existence of $C_0(X)$-valued invariant functionals, playing the role of the Haar measure (Def.\ref{def_is}, Prop.\ref{prop_haar}). 

In Sec.\ref{cx_system}, we introduce the notion of {\em fibred action} of a group bundle $\mG \ra X$ on a $C_0(X)$-algebra (Def.\ref{fibered_act}); we establish some basic properties, and make use of invariant means arising from the invariant $C_0(X)$-valued functionals (Sec.\ref{ss_tp_im}). Usual group actions on \sC algebras are recovered as fibred actions by trivial group bundles (Cor.\ref{cor_uga}).
In Sec.\ref{s_afs}, we discuss the case of fibred actions on Abelian \sC algebras.

In Sec.\ref{dyna}, we provide a generalization of \cite[Thm.1]{DR86} to the non-trivial centre case: every fibred \sC dynamical system $( \mB , \mG )$ is isomorphic to an induced bundle with base space the spectrum $\Omega$ of $\mC$ (Prop.\ref{prop_isob}); moreover, if $\mA' \cap \mB$ is generated by $\mC$ and $\mZ$, then $\mF$ satisfies the analogue of (\ref{eq_rc}) with a group bundle $\mH$ playing the role of $H$ (Thm.\ref{fix_point_3}). Existence of a section $s : X \hra \Omega$ is our main assumption for the above results; this characterizes $\Omega$ as a {\em homogeneous bundle} in the sense of Sec.\ref{cosets}.
We provide a class of examples in the context of the Cuntz algebra, from which it is evident that existence and unicity of $( \mF , \mH )$ are not ensured (Sec.\ref{ss_cuntz}, Sec.\ref{sec_nosec}). 

In Sec.\ref{sec_rc}, we apply our results to Hilbert \sC systems and Doplicher-Roberts endomorphisms associated with superselection sectors. Let $\mA$ be a \sC algebra with centre $\mZ$, and $\rho$ an endomorphism satisfying the special conjugate property in the sense of \cite[\S 2]{DR90}, \cite[\S 4]{DR89A}. Then, we prove that fibred Hilbert \sC systems $(\mF , \mG)$ satisfying $\mF^\mG = \mA$, $\mA' \cap \mF = \mZ$ are in one-to-one correspondence with the set of sections of a suitable homogeneous bundle (Thm.\ref{thm_dr}). In the particular case of endomorphisms studied by Baumg\"artel and Lled\'o, we are also able to prove existence and unicity of $(\mF , \mG)$ (Thm.\ref{thm_bl}), in accord with \cite[Thm.4.13]{BL03}.

%
%

\

{\bf Keywords.} Standard notions about {\em topological} and {\em group bundles} can be found in \cite{Hus}. Let $X$ be a locally compact Hausdorff space; a $C_0(X)$-{\em algebra} is a \sC algebra $\mA$ equipped with a nondegenerate morphism $C_0(X) \ra ZM(\mA)$, where $ZM(\mA)$ is the centre of the multiplier algebra $M(\mA)$ (\cite[\S 2]{Kas88}); in the sequel, we will identify the elements of $C_0(X)$ with their images in $M(\mA)$. $C_0(X)$-algebras correspond to upper-semicontinuous bundles of \sC algebras (\cite[Thm.2.3]{Nil96}), thus generalize the classical notion of {\em continuous field (bundle)} of \sC algebras (\cite[\S 10]{Dix}). If $x \in X$, the {\em fibre} $\mA_x$ of $\mA$ over $x$ is defined as the quotient of $\mA$ by the ideal $\ker x \mA := \left\{ f a \ , \ f \in C_0(X) , f(x) = 0 \ , \ a \in \mA \right\}$. Thus, for every $x \in X$ there is an evaluation epimorphism $\pi_x : \mA \ra \mA_x$, in such a way that $\left\| a \right\| = \sup_x \left\| \pi_x (a)  \right\|$, $a \in \mA$. If $U \subseteq X$ is an open set, then $C_U(X)$ $:=$ $\left\{ f \in C_0(X) : f |_{X-U} = 0  \right\}$ is an ideal of $C(X)$; we define the {\em restriction} $\mA_U$ $:=$ ${\mathrm{closed \ span}}$  $\left\{ fa , f \in C_U(X) , a \in \mA \right\}$, which is a closed ideal of $\mA$. Note that $\mA_U$ is a $C_U(X)$-algebra in the natural way. If $W \subseteq X$ is closed, then the restriction of $\mA$ over $W$ is defined by the epimorphism $\pi_W : \mA \ra \mA_W :=$ $\mA / \mA_{X-W}$.
A $C_0(X)${\em -morphism} from a $C_0(X)$-algebra $\mA$ into a $C_0(X)$-algebra $\mB$ is a \sC algebra morphism $\eta : \mA \ra \mB$ such that $\eta (fa) = f \eta (a)$, $f \in C_0(X)$, $a \in \mA$. Tensor products in the setting of $C_0(X)$-algebras are defined in \cite[\S 3.2]{Bla96}, and are denoted by $\otimes_X$.

\section{Bundles, section groups.}
\label{bundles}

A {\bf bundle} is a surjective, continuous map of locally compact Hausdorff spaces $p : Y \ra X$. We denote by $Y_x := p^{-1}(x)$ the fibre of $Y$ over $x \in X$.

If $U \subseteq X$ is an open set, then $Y_U := p^{-1}(U) \subseteq Y$ defines itself a bundle $p |_{Y_U} : Y_U \ra U$, called the {\em restriction of $Y$ on $U$}. If $p' : Y' \ra X$ is a bundle, a {\em bundle morphism} is a continuous map $F : Y \ra Y'$ such that $p' \circ F = p$; if $F$ is also a homeomorphism, then we say that $Y$, $Y'$ are isomorphic. The {\em fibred product} of bundles $p : Y \ra X$, $p' : Y' \ra X$ is defined as $Y \times_X Y' :=$ $\left\{ (y,y') \in Y \times Y' \ : \right.$ $\left. p(y) = p'(y') \right\}$; note that $Y \times_X Y'$ has a natural structure of a bundle over $X$.
Let $X$, $Y_0$ be locally compact Hausdorff spaces; the natural projection $p : X \times Y_0 \ra X$ defines in a natural way a bundle. A bundle isomorphic to $X \times Y_0$ is called {\em trivial}. 
More in general, a bundle is {\em locally trivial} if for every $x \in X$ there is a neighborhood $U \ni x$ with a bundle isomorphism $\alpha_U : p^{-1} (U) \ra U \times Y_0$, where $Y_0$ is a fixed locally compact Hausdorff space. 

In our terminology, a bundle $p : Y \ra X$ is the topological counterpart of the commutative $C_0(X)$-algebra $C_0(Y)$; the associated structure morphism is $i_p : C_0(X) \ra M(C_0(Y))$, $i_p(f) :=$ $f \circ p$, and the fibres are $C_0(Y)_x =$ $C_0(Y_x)$, $x \in X$; if $z \in C_0(Y)$, the image of $z$ w.r.t. the epimorphism $\pi_x : C_0(Y) \ra C_0(Y_x)$ is given by the restriction $z_x := z |_{Y_x}$.

\begin{defn}
\label{def_cfs}
Let $p : Y \ra X$ be a bundle. A {\bf positive} $C_0(X)${\bf -functional} is a linear, positive $C_0(X)$-module map $\varphi : C_0(Y) \ra C_0(X)$, i.e. $\varphi (fz) = f \varphi (z)$, $\varphi (z') \geq 0$ for every $z \in C_0(Y)$, $f \in C_0(X)$, $0 \leq z' \in C_0(Y)$.
\end{defn}

\begin{lem}
\label{lem_cfs}
Let $p : Y \ra X$ be a bundle, $\varphi : C_0(Y) \ra C_0(X)$ a positive $C_0(X)$-functional. Then, for every $x \in X$ there is a unique positive functional $\varphi_x : C_0(Y_x) \ra \bC$ such that $[\varphi (z)] \ (x) \  =  \  \varphi_x (z_x)$, $z \in C_0(Y)$.
\end{lem}

\begin{proof}
We prove that the map
$\varphi_x  : C_0(Y_x) \ra \bC$, 
$\varphi_x (z_x) :=$ $( \varphi (z) ) \ (x)$, $x \in X$,
is well defined. Let $z' \in C_0(Y)$ such that $z'_x = z_x$, i.e. $z - z' \in \ker \pi_x$, and $C_x(X) \subset C_0(X)$ be the ideal of functions vanishing on $x$; since $\ker \pi_x$ is a nondegenerate Banach $C_x(X)$-bimodule, by \cite[Prop.1.8]{Bla96} there is $f \in C_x(X)$ such that $z - z' = f w$, $w \in \ker \pi_x$. Thus
$\varphi_x (z_x) - \varphi_x (z'_x) =$ 
$(f \varphi (w)) \ (x)  =$
$f(x) \ (\varphi(w)) \ (x) =$
$0$, 
and $\varphi_x$ is well-defined.
\end{proof}

Let $id_X$ denote the identity map over $X$. A (not necessarily continuous) map $s : X \ra Y$ such that $p \circ s = id_X$ is called a {\em selection} of $Y$. In particular, a {\em section} is a {\em continuous} selection $s : X \ra Y$, $p \circ s = id_X$. From the \sC algebraic point of view, sections correspond to $C_0(X)$-epimorphisms
\[  
\phi_s : C_0(Y) \ra C_0(X) \ \ , \ \ z \mapsto \phi_s(z) := z \circ s \ \ .
\]
\noindent In general, existence of a section is not ensured.
%
%
On the other side, a frequently verified property is existence of {\em local sections}, i.e. maps $s_U : U \hra Y$, $p \circ s_U = id_U$, where $U \subset X$ is an open set. We say that $Y$ {\em has local sections} if for every $x \in X$ there is a neighborhood $U \ni x$ with a section $s_U : U \hra Y$. Locally trivial bundles have local sections.

We denote by $\mSY$ the set of sections of $Y$. Now, $\mSY$ is endowed with the "uniform convergence" weak topology such that for every $z \in C_0(Y)$ the map $\left\{ s \mapsto z \circ s \right\}$ is norm continuous. Of course, if $X$ reduces to a single point we get the usual Gel'fand topology over $Y \equiv S_{\left\{ x \right\}} (Y)$. A closed subset $S \subseteq \mSY$ is said to be {\em total} if for every $y \in Y$ there is a section $s \in S$ with $y = s \circ p(y)$. In this case, we say that $Y$ is {\em full}.

A {\em group bundle} is a bundle $p : \mG \ra X$ such that every fibre $G_x := p^{-1} (x)$, $x \in X$, is a locally compact group w.r.t. the topology induced by $\mG$. 

For every $x \in X$, we denote by $e_x$ the identity of $G_x$, and consider the {\em identity selection}
\begin{equation}
\label{def_ids}
e := \left\{ X \ni x \mapsto e_x \in \mG \right\} \ . 
\end{equation}
\noindent  Moreover for every $g , g' \in S_X(\mG)$, the following selections are defined:
\begin{equation}
\label{def_pro}
\left\{
\begin{array}{ll}
gg' \ : \ gg' (x)  := g(x) g'(x)        \\
g^{-1} \ : \ g^{-1} (x) := g(x)^{-1} \ \ , \ x \in X
\end{array}
\right.
\end{equation}

\begin{defn}
\label{def_sg}
Let $p : \mG \ra X$ be a group bundle. A {\bf section group} for $\mG$ is a total subset $G \subseteq \mSG$, which is also a topological group w.r.t. the structure (\ref{def_ids},\ref{def_pro}).
\end{defn}

\begin{ex} 
\label{ex_tb}
Let $G_0$ be a locally compact group and $\mG := X \times G_0$; then, $S_X(\mG)$ can be identified with the group $C(X,G_0)$ of continuous maps from $X$ into $G_0$. It is clear that $S_X(\mG)$ is a section group, so that $\mG$ is full. Let $G \subset \mSG$ be the subgroup of constant maps from $X$ into $G_0$; then, $G$ is a section group, isomorphic to $G_0$.
\end{ex}

%
%
%
%
%
%

\subsection{Homogeneous bundles.}
\label{cosets}

Let $p : \mG \ra X$ be a group bundle, and $\mH \subseteq \mG$ a group subbundle. Then, every fibre $H_x := \mH \cap p^{-1}(x)$, $x \in X$, is a subgroup of $G_x := p^{-1}(x)$. We denote by $\HG$ the quotient space defined by the equivalence relation induced by $\mH$ on $\mG$. For every $y \in \mG$, we denote by $y_H :=$ $\left\{ hy : h \in H_{p(y)}  \right\}$ the associated element of $\HG$. Now, $\HG$ is a locally compact Hausdorff space, endowed with the surjective map
\begin{equation}
\label{def_ph}
p_H : \mG \ra \HG 
\ \ , \ \ 
p_H (y) := y_H 
\ \ .
\end{equation}
\noindent Moreover, $\HG$ has a natural bundle structure $q  : \HG \ra X$, $q (y_H) :=$ $p(y)$, with fibre the space $p_H^{-1} (x) :=$ $H_x \backslash G_x$ of left $H_x$-cosets in $G_x$. Thus, $p_H$ is a bundle morphism, i.e. $q \circ p_H = p$. If $\mG$ is full then the same is true for $\HG$, in fact the set $G_H :=$ $\left\{ g_H := p_H \circ g , g \in \mSG \right\}$ $\subset S_X(\HG)$ is total for $\HG$. 
%
%
To be more concise, we define $\Omega := \HG$, so that we have a bundle $q : \Omega \ra X$. The map (\ref{def_ph}) defines itself a bundle 
$p_\Omega : \mG \ra \Omega$,
that we are going to describe. Let us consider the restriction $p_{\Omega,x} : G_x \ra \Omega_x$, $x \in X$. Then, we may regard $\Omega_x$ as the coset space $H_x \backslash G_x$, so that $p_{\Omega,x}$ can be interpreted as the natural projection from $G_x$ onto $H_x \backslash G_x$, induced by the quotient topology. We have the following commutative diagram
\begin{equation}
\label{dia_hg}
\xymatrix{
   \mH
   \ar[r]^-{\subseteq}
   \ar[dr]_-{p}
 & \mG
   \ar @{>>}[r]^-{p_\Omega}
   \ar[d]_-{p}     
 & \Omega
   \ar[dl]^-{q}    
\\ {}
 & X
  {}
 & 
    }
\end{equation}
\noindent from which it follows that $p_\Omega : \mG \ra \Omega$ has fibres $H_{ p(\omega) } :=$ $G_{p(\omega)} \cap \mH$, $\omega \in \Omega$. Note that if $X$ reduces to a single point, then (\ref{dia_hg}) is simply a quotient map of the type $K \hra$ $K \ra$ $K \backslash L$, where $K \subseteq L$ is an inclusion of compact groups.

%
%

\begin{lem}
\label{lem_lm}
Let $p : \mG \ra X$ be a compact group bundle, $\mH \subseteq \mG$ a compact group subbundle. If $p_\Omega : \mG \ra \Omega$, $\Omega := \HG$, has local sections, then there is a finite open cover $\left\{ W_k \right\}$ of $\Omega$ with bundle morphisms
\begin{equation}
\label{def_lm}
\delta_k : \mG_k := p_\Omega^{-1}(W_k) \ra \mH_k := \mH \cap \mG_k
\end{equation}
\noindent such that 
\begin{equation}
\label{eq_lm}
\delta_k (hy) = h \delta_k(y)
\ , \
h \in \mH_k , y \in \mG_k  \ \ .
\end{equation}
\end{lem}

\begin{proof}
For every $\omega \in \Omega$, there exists a neighborhood $W \subseteq \Omega$, $\omega \in W$, with a section $s_W : W \hra \mG$, $p_\Omega \circ s_W = id_W$. Since $\Omega$ is compact, we obtain a finite open subcover $\left\{ W_k \right\}$, with sections $s_k : W_k \hra \mG$. The morphisms $\delta_k$ are defined by $\delta_k (y) :=$ $y \cdot (s_k \circ p_\Omega (y))^{-1}$.
\end{proof}

\begin{lem}
\label{lem_ls}
Let $p : \mG \ra X$ be a locally trivial bundle with fibre a compact Lie group $L$, $\mH \subseteq \mG$ a compact group subbundle. Then, the bundle $p_\Omega : \mG \ra \Omega$, $\Omega := \HG$ has local sections. 
\end{lem}

\begin{proof}
Let $q : \Omega \ra X$ denote the natural projection, and $\omega \in \Omega$; then, there is a neighborohood $W$ of $\omega$ such that $U := q (W)$ trivializes $\mG$, so that we may identify $\mG_U :=$ $p^{-1}(U)$ with the trivial bundle $U \times L$. 
The previous remarks imply that $\mH_U :=$ $\mG_U \cap \mH$ is "subtrivial", in the following sense: there is a family $\left\{ U_i \subseteq U \right\}_{i \in I}$ such that $U_i \cap U_j =$ $\emptyset$, $\sqcup_i U_i =$ $U$, and $\mH_U$ is the disjoint union $\sqcup_i U_i \times H_i$, where each $H_i$ is a closed subgroup of $L$. This also implies $q^{-1} (U) =$ $\sqcup_i W_i$, $W_i =$ $U_i \times ( H_i \backslash L )$, and $\omega \in \sqcup_i W_i$.
We now proceed with the proof of the existence of local sections. 
By construction, there is a set of indexes $J \subseteq I$ such that $\omega \in$ $\cap_{i \in J} \ovl W_i \simeq$ $\cap_{i \in J} \ovl U_i \times (H_i \backslash L)$, where $\ovl W_i$ (resp. $\ovl U_i$) denotes the closure of $W_i$ (resp. $U_i$).
Let now $i,j \in J$ such that $U_i \cap \ovl U_j$ contains some element $x$. We consider a sequence $\left\{ ( x_n , y_n ) \right\} \subset$ $U_j \times H_j$ converging to $( x , y ) \in$ $\sqcup_i U_i \times H_i$. Since $H_j$ is closed, we find $y \in H_j$; moreover, since $x \in U_i \cap \ovl U_j$, we find $y \in H_i$. Since we may pick an arbitrary converging sequence $\left\{ y_n \right\} \subseteq$ $H_j$, we find that $H_j \subseteq H_i$. By applying the same procedure to every pair of elements of $J$, we conclude that the set $\left\{ H_i \right\}_{i \in J}$ is totally ordered w.r.t. the inclusion, thus there is an index $i_0$ such that $H_i \subseteq H_{i_0}$ for every $i \in J$. In order for more concise notations, we define $H_0 := H_{i_0}$, $U_0 :=$ $U_{i_0}$, $W_0 :=$ $U_0 \times H_0 \backslash L$.
In this way, we have inclusions $U \times (H_0 \backslash L) \subseteq$ $\sqcup_i U_i \times (H_i \backslash L)$; moreover, we may pick $V_0 \subseteq H_0 \backslash L$ with a section $s : V_0 \hra$ $L$ (\cite[Thm.3.4.3]{Hir}). In this way, we obtain a section
$id_U \times s \in$ $S_{U \times V_0} (U \times L)$. 
Since we may identify $U \times V_0$ with a neighborhood of $\omega$, say $W_0$, we conclude that there exists a section from $W_0$ into $U \times L$.
\end{proof}

\section{Actions on $C(X)$-algebras.}
\label{cx_system}

In the present section, we introduce the notion of fibred action of group bundle on a unital $C(X)$-algebra. For this purpose, we briefly recall some results on the categorical equivalence between $C_0(X)$-algebras and a class of topological objects called \sC bundles. In order to simplify the exposition, we consider the case in which $X$ is compact and the given $C(X)$-algebra unital.

\subsection{\sC bundles.}

Let $X$ be a compact Hausdorff space. A {\em \sC bundle} is a Hausdorff space $\Sigma$ endowed with a surjective, open, continuous map $Q : \Sigma \ra X$ such that every fibre $\Sigma_x := Q^{-1}(x)$, $x \in X$, is homeomorphic to a \sC algebra. We assume that $\Sigma$ is full, i.e. for every $\sigma \in \Sigma$ there is a section $a \in S_X(\Sigma)$ such that $a \circ Q (\sigma) = \sigma$. 
Let $Q' : \Sigma' \ra X$ be a \sC bundle. A {\em morphism} from $\Sigma$ into $\Sigma'$ is a continuous map $\phi : \Sigma \ra \Sigma'$ such that: (1) $Q' \circ \phi = Q$; this implies that $\phi (\Sigma_x) \subseteq \Sigma'_x$, $x \in X$; (2) $\phi_x :=$ $\phi |_{\Sigma_x} :$ $\Sigma_x \ra \Sigma'_x$ is a \sC algebra morphism for every $x \in X$. The map $\phi$ is said to be an isomorphism (resp. monomorphism, epimorphism) if every $\phi_x$, $x \in X$, is a \sC algebra isomorphism (resp. monomorphism, epimorphism).

Now, every $a \in S_X (\Sigma)$ defines a family $\left\{ a(x) \in \Sigma_x \right\}$; thus, the set of sections of $\Sigma$ is endowed with a natural structure of $C(X)$-algebra w.r.t. the pointwise-defined *-algebra operations, sup-norm, and pointwise multiplication by elements of $C(X)$. 
On the other side, let $\mA$ be a $C(X)$-algebra with evaluation epimorphisms $\pi_x : \mA \ra \mA_x$, $x \in X$. Then, the set $\wa \mA :=$ $\bigsqcup_{x \in X} \mA_x$, endowed with the natural surjective map $Q : \wa \mA \ra X$, becomes a \sC bundle if equipped with the basis
\begin{equation}
\label{def_tube}
T_{U,a,\eps} :=
\left\{ 
\sigma \in \wa \mA 
\ : \ 
Q(\sigma) \in U
\ {\mathrm{and}} \
\left\| \sigma - \pi_{Q(\sigma)} (a) \right\| < \eps
\right\} \ ,
\end{equation}
where $U \subseteq X$ is open, $a \in \mA$, $\eps > 0$. Some examples follow: the \sC bundle associated with $\mA := C(X) \otimes \mA_0$ is clearly $\wa \mA = X \times \mA_0$, with $Q(x,a_0) := x$, $x \in X$, $a_0 \in \mA_0$.; if $\mA := C(Y)$ is Abelian, then there is a surjective map $q : Y \ra X$, and the \sC bundle $\wa \mA \ra X$ has fibres $\mA_x = C(q^{-1}(x))$, $x \in X$.

The map $\left\{ \right. \mA \mapsto \wa \mA \left. \right\}$ defines a functor: if $\eta : \mA \ra \mA'$ is a $C(X)$-algebra morphism, then the associated morphism $\wa \eta : \wa \mA \ra \wa{\mA'}$ is defined in such a way that $\wa \eta \circ \pi_x (a) =$ $\pi'_x \circ \eta (a)$, $a \in \mA$, $x \in X$, where $\pi'_x : \mA' \ra \mA'_x$ denotes the evaluation epimorphism.
For the proof of the following theorem, see \cite[Thm.5.13]{Gie82}; in order to understand the terminology used in the above reference, we recommend the reader to \cite[\S 1 , \S 5]{Gie82}.

\begin{thm}
\label{thm_ss}
Let $X$ be a compact Hausdorff space. The functor $\left\{ \right. \mA \mapsto \wa \mA \left. \right\}$ induces an equivalence from the category of $C(X)$-algebras into the category of \sC bundles with base space $X$.
\end{thm}

For every $C(X)$-morphism $\eta : \mA \ra \mA'$, it turns out that $\ker \eta$ is a $C(X)$-algebra. We denote by $\ker \wa \eta \ra X$ the associated \sC bundle. It is clear that there is an inclusion $\ker \wa \eta \subseteq \wa \mA$ such that $(\ker \wa \eta)_x \subseteq$ $\ker \eta_x$, $x \in X$.

We conclude the present section with a remark on $C(X)$-algebras and \sC bundles. Let $\mA$ be a unital $C(X)$-algebra with fibre epimorphisms $\pi_x : \mA \ra \mA_x$, $x \in X$, and $p : Y \ra X$ a compact bundle, so that $C(Y)$ is a unital $C(X)$-algebra. Let $P : \wa \mA \ra X$ be the \sC bundle associated with $\mA$. We consider the set $S_X (Y,\wa \mA)$ of continuous maps $F : Y \ra \wa \mA$ such that
$P \circ F = p$,
$\left\| F \right\| :=$
$\sup_y \left\| F(y) \right\| <$
$\infty$.
$S_X(Y,\wa \mA)$ is endowed with a natural structure of \sC algebra w.r.t. the pointwise-defined operations, and it is easily verified that there is a natural isomorphism
\begin{equation}
\label{eq_pb}
S : C(Y) \otimes_X \mA  \ra  S_X(Y,\wa \mA)
\ \ , \ \ 
[S(z \otimes a)] \ (y) :=  z(y) \ \pi_{p(y)}(a) 
\ \ . 
\end{equation}

\subsection{Group Bundles acting on $C(X)$-algebras.}

Let $G$ be a locally compact group; as remarked in \cite[\S 4]{Nil96}, every action $\alpha : G \ra \aut_X \mA$ induces a family of actions $\left\{ \alpha^x : G \ra \aut \mA_x \right\}_x$, such that
\begin{equation}
\label{eq_fa}
\pi_x \circ \alpha_g = \alpha^x_g \circ \pi_x 
\ \ , \ \ 
x \in X \ , \ a \in \mA \ , \ g \in G \ .
\end{equation}
\noindent The above elementary remark suggests to introduce the following notion.

\begin{defn}
\label{fibered_act}
Let $X$ be a compact space, $p : \mG \ra X$ a group bundle, $\mA$ a $C(X)$-algebra. We say that $\mG$ acts fiberwise on $\mA$ if there is a family of strongly continuous actions $\alpha :=$ $\left\{ \alpha^x : G_x \ra \right.$ $\aut \mA_x $,$\left. x \in X \right\}$, such that the map
\begin{equation}
\label{eq_axg00}
\alpha : \ \mG \times_X \wa \mA \ra \wa \mA
\ \ , \ \
(y,\sigma) \mapsto \alpha^x_y (\sigma)
\ \ , \ \ 
x := p(y) = Q(\sigma)
\ ,
\end{equation}
\noindent is continuous. In such a case, we say that $\alpha$ is a {\bf fibred action}, and that the triple $(\mA , \mG , \alpha )$ is a {\bf fibred} $C(X)$-system. If $\mG$ is full, then we say that $\alpha$ is full. The {\bf fixed-point algebra w.r.t.} $\alpha$ is the $C(X)$-subalgebra $\mA^\alpha :=$ $\left\{ a \in \mA : \alpha^x_y \circ \pi_x(a) = \pi_x(a) \ \forall x \in X \ , \ \forall y \in G_x \right\}$.
\end{defn}

\begin{prop}
\label{prop_ffa}
Let $( \mA , \mG , \alpha )$ be a full fibred $C(X)$-system. Then, for every section group $G \subseteq \mSG$ there is a strongly continuous action 
\begin{equation}
\label{eq_fib_act}
\alpha^G : G \ra {\bf aut}_X \mA 
\end{equation}
\noindent such that 
\begin{equation}
\label{eq_axg0}
\alpha^x_{g(x)} \circ \pi_x = \pi_x \circ \alpha^G_g
\ \ , \ \ g \in G \ , \ x \in X \ .
\end{equation}
\noindent On the converse, if a strongly continuous action (\ref{eq_fib_act}) endowed with a family $\left\{ \alpha^x : G_x \right.$ $\left. \ra {\bf aut} \mA_x \right\}_x$ satisfying (\ref{eq_axg0}) is given, then there is a unique fibred $C_0(X)$-system $( \mA , \mG , \alpha )$ such that $\alpha^G$ is associated with $\alpha$ in the above sense.
\end{prop}

\begin{proof}
Suppose that a fibred action $\alpha$ is given. If $g \in S_X(\mG)$ and $a \in \mA$, then by continuity of $\alpha$ the map $s_{g,a} : X \ra \wa \mA$, $s_{g,a}(x) := \alpha^x_{g(x)} \circ \pi_x (a)$, belongs to $S_X(\wa \mA)$. Elementary remaks show that the map $\alpha^G_g (a) :=$ $s_{g,a}$ defines the desired automorphic action.
Let now $\alpha^G$ be an action satisfying (\ref{eq_fib_act},\ref{eq_axg0}). Since $\mG$ is full, for every $y \in \mG$ there is $g \in G$ with $g(p(y)) = y$. We define the map $\alpha : \mG \times_X \wa \mA \ra \wa \mA$, $\alpha (g(x),\sigma) := \alpha^x_{g(x)}(\sigma)$, $\sigma \in \mA_x$. Note that (\ref{eq_axg0}) ensures that $\alpha ( y,\sigma )$, $y \in \mG$, $\sigma \in \wa \mA$, is well-defined, in fact it does not depend on the choice of $g \in G$ satisfying $g(x) = y$.
We now verify that $\alpha$ is continuous. For this purpose, let us consider $(y,\sigma) \in$ $\mG \times_X \wa \mA$, with $x :=$ $p(y) =$ $Q(\sigma)$; we pick $a \in \mA$ and $g \in G$ such that $g(x) = y$ and $\pi_x (a) =$ $\sigma$, so that $\pi_x \circ \alpha^G_g (a) =$ $\alpha^x_y(\sigma)$. Moreover, we define $b := \alpha^G_g(a)$, and consider a neighborhood of the type $T_{U,b,\eps}$ for $\alpha^x_y(\sigma)$, with $U \subset X$, $U \ni x$ (see (\ref{def_tube})). 
Since $\alpha^G$ is strongly continuous, there is a neighborhood $V \ni y$ with the following property: for every $g' \in G$ such that $g(U) \subseteq V$, and for every cutoff $\lambda \in$ $C_W(X)$, $U \subset W$, $\lambda |_U \equiv 1$, it turns out that 
$\left\|  \lambda b  - \lambda \alpha^G_{g'} (a)  \right\| <$ $\eps / 2$,
so that for every $x' \in U$ we find
$\left\|  \pi_{x'} (b)  - \alpha^{x'}_{y'} \circ \pi_{x'} (a)  \right\| <$ $\eps / 2$,
$y' :=$ $g'(x') \in V$.
Thus, for every $( y' , \sigma' ) \in$ $V \times_X T_{U,a,\eps / 2}$, it turns out, with $x' :=$ $p(y') =$ $Q(\sigma') \in U$,
\[
\left\| \pi_{x'} (b)  - \alpha^{x'}_{y'} (\sigma')  \right\| \leq
\left\| \pi_{x'} (b)  - \alpha^{x'}_{y'} \circ \pi_{x'} (a)  \right\| +
\left\| \pi_{x'} (a) - \sigma'  \right\| \leq
\eps / 2 +
\eps / 2
\ \ .
\]
\noindent Thus, $\alpha^{x'}_{y'} (\sigma') \in$ $T_{U,b,\eps}$, and we conclude that $\alpha$ is continuous.
\end{proof}

\begin{cor} 
\label{cor_uga}
Let $\mA$ be a unital \sC algebra, $G$ a locally compact group, and $\alpha : G \ra \aut \mA$ a strongly continuous action. We define $C(X) := (\mA \cap \mA')^\alpha$ (so that $\mA$ is a $C(X)$-algebra) and $\mG := X \times G$. There is a unique fibred action $\wt \alpha : \mG \times_X \wa \mA \ra \wa \mA$ such that $\alpha$ is the action associated with $\wt \alpha$ in the sense of the previous proposition.
\end{cor}

\begin{proof}
$G$ appears as a section group of constant sections for $\mG := X \times G$, so that (\ref{eq_fa}) can be regarded as a special case of (\ref{eq_axg0}). The proof now follows from Prop.\ref{prop_ffa}.
\end{proof}

%
%

We introduce a natural notion of equivariance. If $( \mA , \mG , \alpha )$, $( \mB , \mG , \beta )$ are fibred $C(X)$-systems and $\eta : \mA \ra \mB$ is a $C(X)$-morphism, we say that $\eta$ is $\mG$-{\em equivariant} if $\wa \eta \circ \alpha ( y , \sigma ) =$ $\beta ( y , \wa \eta (\sigma) )$, $y \in \mG$, $\sigma \in \wa \mA$. In this case, we use the notation $\eta :$ $(\mA,\alpha) \ra$ $(\mB,\beta)$. The proof of the following Lemma is trivial, therefore we omit it.

\begin{lem}
\label{lem_qf}
Let $p : \mG \ra X$ a group bundle, $( \mA , \mG , \alpha )$ a fibred $C(X)$-system, and $\eta : \mA \ra \mF$ a $C(X)$-algebra epimorphism with associated family $\left\{ \eta_x : \mA_x \ra \right.$ $\left. \mF_x \right\}_{x \in X}$. If $\mH \subseteq \mG$ is a group subbundle such that
\begin{equation}
\label{def_st}
\alpha (y,\sigma) \in \ker \wa \eta
\ \ , \ \
y \in \mH \ , \ \sigma \in \ker \wa \eta \ ,
\end{equation}
\noindent then there is a unique fibred action $\beta : \mH \times_X \wa \mF \ra \wa \mF$ such that $\wa \eta \circ \alpha ( y , \sigma ) =$ $\beta ( y , \wa \eta (\sigma) )$, $y \in \mH$, $\sigma \in \wa \mA$.
\end{lem}

%
%

\begin{cor}
\label{cor_fa_q}
Let $\mA$ be a $C(X)$-algebra, $G$ a compact group, and $\alpha : G \ra \aut_X \mA$ a strongly continuous action. If $\eta : \mA \ra \mF$ is a $C(X)$-epimorphism, then every element of the stabilizer $H :=$ $\left\{ h \in G : \alpha_h (\ker \eta) = \ker \eta \right\}$ defines a section of the group bundle $\mH \subseteq$ $X \times G$ of elements satisfying (\ref{def_st}) w.r.t. the induced fibred action on $\mA$.
\end{cor}

\begin{rem}
\label{rem_fa_q}
According to the previous corollary, there is a fibred action $( \mF , \mH ,$ $ \beta )$, anyway also a strongly continuous action $\beta_H : H \ra \aut_X \mF$ is defined. Note that in general $H$ is not a section group for $\mH$; for example, if we consider $X = [0,1]$, $\mH =$ $\left\{ (x,h) \in X \times G : h = 1 \ \mathrm{for} \ x \geq 1/2 \right\}$, then we obtain $H = \left\{ 1 \right\}$. Thus, we may have the "extreme" situation in which $H$ reduces to the trivial group, whilst $\mH$ is non-trivial (this happens in the case in which $S_X(\mH)$ reduces to the identity selection). Roughly speaking, fibred actions are better behaved than usual strongly continuous actions w.r.t. $C(X)$-epimorphisms.
\end{rem}

%
%
%

%
%

\begin{ex}
\label{ex_ta}
Let $p : \mG \ra X$ be a compact group bundle. We denote by $C(\mG)$ the $C(X)$-algebra of continuous functions of $\mG$, and by $Q : \wa \mA_\mG \ra$ $X$ the associated \sC bundle with fibres $\mA_{\mG,x} = C(G_x)$, $x \in X$.
It is clear that there are actions $\lambda^x : G_x \ra$ $\aut C(G_x)$, $\lambda^x_h \xi (y) :=$ $\xi (h^{-1}y)$, $\rho^x : G_x \ra$ $\aut C(G_x)$, $\rho^x_h \xi (y) :=$$\xi (yh)$, $\xi \in C(G_x)$, $h,y \in G_x$.
Thus, there are fibred actions
\begin{equation}
\label{aut_g}
\left\{
\begin{array}{ll}
\lambda : \mG \times \wa \mA_\mG \ra \wa \mA_\mG 
\ \ , \ \  
\lambda ( y , \xi ) := \lambda^{p(y)}_y (\xi)
\\
\rho : \mG \times \wa \mA_\mG \ra \wa \mA_\mG 
\ \ , \ \  
\rho ( y , \xi ) := \rho^{p(y)}_y (\xi)
\end{array}
\right.
\end{equation}
\noindent We call $\lambda$, $\rho$ respectively the left and right
translation actions over $\mG$. 
\end{ex}

%
%

%
%

%
%
%

Let $p : \mG \ra X$ be a group bundle, and $( \mA , \mG , \alpha )$ a fibred $C(X)$-system. We denote by $X'$ the spectrum of $(\mA' \cap \mA)^\alpha$, and identify $C(X')$ with $(\mA' \cap \mA)^\alpha$. The obvious inclusion $C(X) \hra C(X')$ induces a surjective map $q : X' \ra X$. Now, it is clear that $\mA$ is a $C(X')$-algebra with associated \sC bundle $Q' : \wa \mB \ra X'$, and there is a $C(X')$-isomorphism $\tau : \mA \ra $ $\mB := S_{X'}(\wa \mB)$. At the level of \sC bundles, it is easy to recognize that $\tau$ induces an epimorphism
$\wt \tau : X' \times_X \wa \mA \ra \wa \mB$.
%
%
%
%
Let $q_* \mG :=$ $X' \times_X \mG$; then, the following map is well-defined:
\begin{equation}
\label{rem_ce}
\beta : q_* \mG \times_{X'} \wa \mB \ra \wa \mB 
\ \ , \ \ 
\beta ( (x',y) \ , \ \wt \tau ( x' , \sigma) )  
\ := \ 
\wt \tau ( x' , \alpha ( y , \sigma) )
\ \ . 
\end{equation}
\noindent By construction $( \mB , q_* \mG , \beta )$ is a fibred $C(X')$-system, and $(\mB' \cap \mB)^\beta = C(X')$. We conclude that fibred $C(X)$-systems $( \mA , \mG , \alpha )$ can be always "rescaled" in such a way that $C(X)$ coincides with $(\mA' \cap \mA)^\alpha$.

\section{Invariant means.}
\label{ss_tp_im}

We recall the reader to the notation used in Lemma \ref{lem_cfs} and Ex.\ref{ex_ta}.

\begin{defn}
\label{def_is}
Let $p : \mG \ra X$ be a group bundle. A positive $C_0(X)$-functional $\varphi : C_0 (\mG) \ra C_0 (X)$ is said to be {\em left (resp. right) invariant} if
\begin{equation}
\label{eq_is}
\varphi_x \circ \lambda ( y , z_x )   =  \varphi_x (z_x) 
\ \ {\mathrm{(resp.}} \ \ 
\varphi_x \circ \rho ( y , z_x ) = \varphi_x (z_x) \ ) \ ,
\end{equation}
\noindent $y \in \mG$, $x := p(y)$, $z \in C_0(\mG)$. If $\varphi$ is left and right invariant, then $\varphi$ is said to be {\bf invariant}.
\end{defn}

In the next lemma, we establish a unicity result for invariant functionals. We denote by $C_b(X)$ (= $M(C_0(X))$) the \sC algebra of bounded continuous functions on $X$, and by $C_b^+(X)$ the space of positive functions.

\begin{lem}
\label{lem_is}
Let $p : \mG \ra X$ be a group bundle, $\varphi$ a left (resp. right) invariant positive $C_0(X)$-functional. Let $\mu_x : C_0(G_x) \ra \bC$ denote a left (resp. right) Haar measure on $G_x$, $x \in X$. Then, $\varphi$ is unique up to a multiplicative factor $\tau \in C_b(X)^+$.
\end{lem}

\begin{proof}
Let $\left\{ \varphi_x \in C_0(G_x)^* \right\}_x$ be the family of states defined as in Lemma \ref{lem_cfs}. The equality (\ref{eq_is}) implies that $\varphi_x$ is a left (resp. right) $G_x$-invariant state on $C_0(G_x)$; the left (resp. right) $G_x$-invariance of $\varphi_x$ implies that $\varphi_x = \chi_x \mu_x$ for some $\chi_x \in \bR^+$. This also implies that if $\varphi' : C_0(\mG) \ra C_0(X)$ is a left (resp. right) invariant positive $C_0(X)$-functional then we find 
$\chi'_x \cdot (\varphi (z)) \ (x) =$ $\chi_x \cdot (\varphi'(z)) (x)$,
$z \in C_0(\mG)$. The previous equality implies that $\varphi (z)$, $\varphi'(z)$ $\in C_0(X)$ have the same support for every $z \in \mZ$. We define the map $\tau : X \ra \bC$, $\tau (x) := \chi_x^{-1} \chi'_x$, $x \in X$, and prove that it is continuous; for this purpose, we consider an approximate unit $\left\{ \lambda_i \right\}_{i} \subset C_0(\mG)$, and note that $\tau$ is a limit in the strict topology $\tau =$ $\lim_i \varphi'(\lambda_i) \varphi(\lambda_i)^{-1}$;
since the net in the r.h.s. of the previous equality is contained in $C_0(X)$, we conclude that $\tau \in C_b(X)$.
\end{proof}

Let $G_0$ be a locally compact group, $K$ a topological group acting continuously by proper homeomorphisms on $G_0$. Then, $K$ acts by automorphisms on $C_0(G_0)$; in the sequel, we will identify elements of $K$ with the corresponding automorphisms of $C_0(G_0)$.

\begin{prop}
\label{prop_haar}
Let $X$ be a locally compact, paracompact Hausdorff space, $p : \mG \ra X$ a locally trivial group bundle with fibre $G_0$ and structure group $K$. Let $\varphi_0 : C_0(G_0) \ra \bC$ denote a left (resp. right) Haar measure such that $\varphi_0 \circ \alpha = \varphi_0$, $\alpha \in K$. Then, there is a left (resp. right) invariant positive $C_0(X)$-functional on $C_0(\mG)$.
\end{prop}

\begin{proof}
Let $\left\{ U_i \right\}_{i \in I}$ be a locally finite, trivializing open cover for $\mG$. Then $C_0(\mG)$ is a locally trivial continuous bundle with local charts $\alpha_i : C_0(p^{-1}(U_i)) \ra$ $C_0(U_i) \otimes C_0(G_0)$, and $K$-cocycle
$\alpha_{ij} :$ $U_i \cap U_j \ra K \subseteq$ $\aut C_0(G_0)$. 
Let $id_i$ be the identity on $C_0(U_i)$; we define $\varphi_i : C_0(p^{-1}(U_i)) \ra C_0(U_i)$,  $\varphi_i := ( id_i \otimes \varphi_0 ) \circ \alpha_i$. It is clear that every $\varphi_i$ is a left (resp. right) invariant $C_0(U_i)$-functional on $C_0(p^{-1}(U_i))$. Now, since $\varphi_0$ is $K$-invariant we find that if $U_i \cap U_j \neq \emptyset$, then
$(id_{ij} \otimes \varphi_0)  \circ  ( id_{ij} \otimes \alpha ) =$ 
$id_{ij} \otimes \varphi_0$, $\alpha \in K$. This implies that $\varphi_i (F) = \varphi_j (F)$, $F \in C_0(p^{-1}(U_i)) \cap C_0(p^{-1}(U_j))$.
Thus, if $\left\{ \lambda_i  \right\}$ is a partition of unity subordinate to $\left\{ U_i \right\}$, then the following left (resp. right) invariant positive $C_0(X)$-functional is defined: $\varphi : C_0(\mG) \ra C_0(X)$, $\varphi (F) := \sum_i \varphi_i (\lambda_i F)$, $F \in C_0(\mG)$.
\end{proof}

It is clear that every trivial group bundle admits an invariant $C_0(X)$-functional. Moreover, by applying the previous result we may also easily prove that the same is true for principal bundles, and for bundles of unitary, vector-bundle morphisms (\cite[I.4.8]{Kar}).

%
%
%
%
%
%

Let $\mA$ be a $C(X)$-algebra with fibre epimorphisms $\pi_x : \mA \ra$ $\mA_x$, $x \in X$, and $p : \mG \ra X$ a compact group bundle. To be concise, we define $\mA^\mG :=$ $\goa$. $\mA^\mG$ has a natural structure of $C(X)$-algebra, with fibres $\mA^\mG_x =$ $C(G_x) \otimes \mA_x$, $x \in X$. We denote by $P : \wa \mA^\mG \ra X$ the associated \sC bundle.
The translation actions defined in Ex.\ref{ex_ta} extend in a natural way on $\mA^\mG$:
\begin{equation}
\label{def_ta}
\left\{
\begin{array}{ll}
\lambda^\mA : \mG \times \wa \mA^\mG \ra \wa \mA^\mG
\ \ , \ \ 
\lambda^{\mA,x}_y F (y') := F (y^{-1}y')
\\
\rho^\mA : \mG \times \wa \mA^\mG \ra \wa \mA^\mG
\ \ , \ \ 
\rho^{\mA,x}_y F (y') = F (y'y)
\end{array}
\right.
\end{equation}
\noindent where $F \in $ $S_X(\mG,\wa \mA) \simeq$ $\mA^\mG$ (see (\ref{eq_pb})), $y , y' \in \mG$, $x := p(y) = p(y')$.
We call (\ref{def_ta}) {\em left (resp. right) translation} $G$-{\bf action} over $\goa$. The following lemma is just the translation in terms of fibred systems of a well-known basic property of \sC algebra actions, thus we omit the proof.

\begin{lem}
\label{lem_mfa}
Let $p : \mG \ra X$ be a compact group bundle, and $( \mA , \mG , \alpha )$ a fibred $C(X)$-system. Then, there is an equivariant $C(X)$-monomorphism
$T  : ( \mA , \alpha ) \ra$ $( \mA^\mG , \rho^\mA )$ such that 
$Ta (y) :=$ $\alpha ( y , \pi_x (a) )$, 
$a \in \mA$, $y \in \mG$, $x := p(y)$.
\end{lem}

Let $p : \mG \ra X$ be a compact group bundle with a right invariant positive $C(X)$-functional $\varphi : C(\mG) \ra C(X)$, and $\mA$ a $C(X)$-algebra. Then, the following positive functional is defined
\[
\varphi_\mA : \mA^\mG \ra \mA \ \ : \ \ 
\varphi_\mA (z \otimes a) := \varphi (z) \ a  \ ,
\]
\noindent in such a way that the equality
$\pi_x \circ \varphi_\mA (F)  =$ $\int_{G_x} F(y) \ d\mu_x (y)$
holds for every $F \in \mA^\mG$, $x := p(y) \in X$; $\mu_x$ denotes a Haar measure of $G_x$, according to Lemma \ref{lem_is}. 
Let $\alpha : \mG \times_X \wa \mA \ra$ $\wa \mA$ be a fibred action. Then, a $\mG$-invariant mean $m_\mG : \mA \ra \mA^\alpha$, $m_\mG (a) :=$ $\varphi_\mA \circ Ta$, is defined, where $Ta \in$ $\mA^\mG$ is defined as in Lemma \ref{lem_mfa}. For every $x \in X$, the following equality holds:
\begin{equation}
\label{eq_int}
\pi_x \circ m_\mG (a) = \int_{G_x} \alpha (  y , \pi_x (a) ) \ d \mu_x (y)
\ \ .
\end{equation}
\noindent We define
$
m_{\mG , x} : \mA_x \ra (\mA_x)^{\alpha^x} 
$
as the l.h.s. of (\ref{eq_int}), so that $m_{\mG,x} \circ \pi_x = \pi_x \circ m_\mG$.

\section{Abelian fibred $C(X)$-systems.}
\label{s_afs}

Structural properties of centers of \sC algebras carrying a group action will be important in the sequel. For this reason, we establish some properties of fibred actions in the Abelian case.

We fix some notations for the rest of the present section. Let $X$ be a compact Hausdorff space, $\mC := C(\Omega)$ a unital, Abelian $C(X)$-algebra, so that there is a surjective map $q : \Omega \ra X$; we define $\Omega_x := q^{-1}(x)$, $x \in X$. In such a way, the $C(X)$-algebra structure of $\mC$ is described by the monomorphism $i_q : C(X) \ra \mC$, $i_q (f) := f \circ q$. We denote the evaluation epimorphisms by $\eta_x : \mC \ra \mC_x \simeq C(\Omega_x)$, $x \in X$.

Let $p : \mG \ra X$ be a compact group bundle. We consider a fibred action $\alpha : \mG \times_X \wa \mC \ra$ $\wa \mC$, such that $\mC^\alpha =$ $i_q (C(X)) \simeq$ $C(X)$. By applying the Gel'fand transform, we obtain an action by bundle automorphisms
\begin{equation}
\label{def_}
\alpha^* : \mG \times_X \Omega \ra \Omega 
\ \ , \ \ 
\alpha^* (y , \omega ) := \omega \circ \alpha^x_y  \ , \ x := q(\omega)
\ .
\end{equation}
\noindent Let $( C(\Omega')  , \mG , \beta )$ be a fibred $C(X)$-system. A (bundle) morphism $F : \Omega \ra \Omega'$ is $\mG$-{\em equivariant} if $F \circ \alpha^* (y,\omega) =$ $\beta^* ( y , F(\omega) )$, $y \in \mG$, $\omega \in \Omega$; in this case, we write $F :$ $(\Omega , \alpha^*) \ra$ $( \Omega' , \beta^* )$. 

In the sequel of the present section, {\em we assume that there exists a left invariant positive $C(X)$-functional} $\varphi : C(\mG) \ra C(X)$.

\begin{lem}
\label{stru_c}
$\mC$ is a continuous bundle of Abelian \sC algebras over $X$; for every $x \in X$, there is a closed group $H_x \subseteq G_x$, unique up to isomorphism, such that the fibre $\mC_x$ is isomorphic to $C(H_x \backslash G_x )$.
\end{lem}

\begin{proof}
Let $\omega \in \Omega$, $x := p (\omega) \in X$. For every $x \in X$ there is an action $\alpha^x : G_x \ra$ $\aut \mC_x$, with $\mC_x \simeq$ $C(\Omega_x)$.
We now prove that $G_x$ acts transitively on $\Omega_x$. Let $f \in C(\Omega_x)$ be $G_x$-invariant; by the Tietze theorem, there is $f' \in C(\Omega)$ with $f'|_{\Omega_x} = f$; by applying the $\mG$-invariant mean $m_G : \mC \ra$ $\mC^\alpha$, we obtain a function ${f_0} := m_\mG (f') \in \mC^\alpha$. By (\ref{eq_int}), we obtain
${f_0} |_{\Omega_x} =$ $m_\mG (f') |_{\Omega_x} =$ $m_{\mG , x} (f' |_{\Omega_x} ) =$
$m_{\mG , x} (f) =$ $f$.
Since ${f_0} \in$ $i_q (C(X)) =$ $\mC^\alpha$, we obtain that ${f_0}$ is a constant function on $\Omega_x$ for every $x \in X$. This proves that $f$ is a constant map; thus $C(\Omega_x)^{\alpha^x} = \bC 1$, i.e. $G_x$ acts transitively on $\Omega_x$. The isomorphism $C(\Omega_x) \simeq C(H_x \backslash G_x)$ is obtained by defining $H_x$ as the stabilizer of some $\omega_x \in \Omega_x$. 
It remains to verify that $\mC$ is a continuous bundle over $X$. For this purpose, we consider the right Hilbert $C(X)$-module $\mM := L^2 (\mC , m_\mG)$, defined by completition of $\mC$ w.r.t. the $C(X)$-valued scalar product $\langle c , c' \rangle :=$ $m_\mG (c^* c')$; by construction, there is a $C(X)$-module map $\mC \ni c \mapsto v_c \in \mM$. 
For every $x \in X$, we consider the state $\omega_x := x \circ m_\mG : \mC \ra \bC$, and the associated GNS representation $\wt \pi_x : \mC \ra L(\mM_x)$. By definition, $\mM_x$ coincides with the fibre of $\mM$ over $x$ in the sense of \cite[Def.14]{Bla96}.
Now, $\mC$ acts by multiplication on $\mM$, so that there is a $C(X)$-morphism $\pi : \mC \ra L(\mM)$, $\pi (c) v_{c'} := v_{cc'}$, $c,c' \in \mC$. Since $m_\mG$ is positive, we find that $\pi$ is injective. 
If $c,c' \in \mC$ and $\eta_x (c - c') = 0$, then there is $f \in C(X)$ such that $f(x) = 0$ and $c - c' = i_q (f) c''$ for some $c'' \in \mC$ (see \cite[Prop.1.8]{Bla96}). This implies $\omega_x ( c - c' ) = f(x) m_\mG (c'') = 0$, i.e. $\wt \pi_x (c - c') = 0$.
We conclude that $\eta_x (c) = \eta_x (c')$ $\Rightarrow$ $\wt \pi_x (c) = \wt \pi_x (c')$; thus, $\wt \pi_x = \pi_x \circ \eta_x$ for some representation $\pi_x : \mC_x \ra L(\mM_x)$.
The previous considerations imply that $\pi$ defines a field $\left\{ \pi_x \right\}_{x \in X}$ of faithful representations in the sense of \cite[Def.2.11]{Bla96}, and \cite[Prop.3.11]{Bla96} implies that $\mC$ is a continuous bundle of \sC algebras over $X$. 
\end{proof}

\begin{lem}
\label{iso_c}
Suppose that there is a section $s \in S_X(\Omega)$. Then, there is a compact group subbundle $\mH \subseteq$ $\mG$, with a bundle isomorphism $\Omega \simeq$ $\HG$ (see Sec.\ref{cosets}).
\end{lem}

\begin{proof}
We define $\mH :=$ $\left\{ h \in \mG : \right.$ $\left. \alpha^* ( \ h , s(p(h)) \ ) = s(p(h)) \right\}$; moreover, we note that for every $y \in \mG$, $h \in \mH$, it turns out $\alpha^* ( \ hy , s(p(hy)) \ ) =$ $\alpha^* ( \ y , s(p(y)) \ )$; thus, the map 
\begin{equation}
\label{def_tau}
\tau : \HG \ra \Omega
\ \ , \ \
\tau (y_H) := \alpha^* ( \ y , s(p(y))  \ )
\end{equation}
is well-defined, and it is trivial to check that it is a bundle isomorphism.
\end{proof}

We conclude that if $( C(\Omega) , \mG , \alpha )$ is a fibred $C(X)$-system endowed with an invariant $C(X)$-functional $\varphi : C(\mG) \ra C(X)$ and a section $s \in S_X(\Omega)$, then there exists a compact group subbundle $\mH_s \subseteq \mG$ with an isomorphism $\HsG \simeq \Omega$. Thus, the considerations of Sec.\ref{cosets} apply, in particular Lemma \ref{lem_lm} and Lemma \ref{lem_ls}. Note that the isomorphism class of $\mH_s$ depends on $s$: by choosing a different section $s' \in S_X(\Omega)$, we may get a group bundle $\mH_{s'}$ not isomorphic to $\mH_s$ (Sec. \ref{sec_nosec}).

\section{Induced \sC bundles.}
\label{dyna}

In the present section, we denote by $p : \mG \ra X$ a compact group bundle. Moreover, we consider a fibred $C(X)$-system $( \mB , \mG , \alpha )$. We fix our notation with $\mA := \mB^\alpha$, $\mZ := \mA \cap \mA'$, $\mC := \mB \cap \mB'$. By (\ref{rem_ce}), we may eventually "rescale" $X$, and assume that $C(X) =$ $\mC^\alpha$.
We denote by $\mC \vee \mZ$ the \sC subalgebra of $\mB$ generated by $\mC$ and $\mZ$. We denote by $\pi_x : \mB \ra \mB_x$, $x \in X$, the epimorphisms associated with $\mB$ as an upper semicontinuous bundle over $X$. Note that $\mA$, $\mZ$, $\mC$ inherit from $\mB$ the structure of $C(X)$-algebra. We denote by $\Omega$ the spectrum of $\mC$, and by $q : \Omega \ra X$ the natural projection induced by the $C(X)$-algebra structure of $\mC$; in particular, we define $\Omega_x :=$ $q^{-1}(x)$, $x \in X$. Since $\mC$ is stable w.r.t. $\alpha$, by applying the Gel'fand transform we obtain a fibred action $\alpha^* : \mG \times \Omega \ra \Omega$, $\alpha^* (y,\omega) :=$ $\omega \circ \alpha^{p(y)}_y$.
We now expose the main results of the present section. Together with Prop.\ref{prop_isob}, they generalize \cite[Thm.1]{DR86} to the nontrivial centre case.

\begin{thm}
\label{fix_point_3}
Suppose that $\mG$ admits a left invariant positive $C(X)$-functional, and consider the fibred $C(X)$-system $( \mB , \mG , \alpha )$. Then, for every section $s \in S_X(\Omega)$ there is a compact group subbundle $\mH_s \subseteq \mG$ and a fibred $C(X)$-system $({\mF_s}, {\mH_s},\beta)$, with a $C(X)$-epimorphism $\eta : \mB \ra \mF_s$. Moreover, the following properties are satisfied:
\begin{enumerate}
\item  $\mC$ is a continuous bundle of \sC algebras over $X$, with fibres 
       $C(H_x \backslash G_x)$, $x \in X$; there is a $C(X)$-algebra 
       isomorphism $\mC \simeq$ $C(\HsG)$;
\item  $\eta$ is injective on $\mA$, and $\mA \simeq$ $\eta(\mA) =$ ${\mF_s}^\beta$;
\item  assume that $\mA' \cap \mB = \mC \vee \mZ$; if $\mG \ra \Omega$ has local sections, 
       or $\mH_s$ is endowed with a right-invariant positive $C(X)$-functional,
       then $\mZ \simeq$ $\eta(\mZ) =$ $\eta(\mA)' \cap {\mF_s}$.
\end{enumerate}
\end{thm}

Some remarks follow.
Existence of a (left) right invariant positive $C(X)$-functional is verified under quite general conditions, in particular when $\mG \ra X$ is trivial (see Prop.\ref{prop_haar}); recall that compact group actions on \sC algebras in the usual sense can be regarded as fibred actions by trivial group bundles (Cor.\ref{cor_uga}).
Existence of the section $s \in S_X(\Omega)$ is the necessary assumption to define the group bundle $\mH_s$, on which our construction is based; fibred $C(X)$-systems such that $\Omega \ra X$ does not admit sections can be constructed also in the case in which $\mG \ra X$ is trivial, as we shall see in Sec.\ref{sec_nosec}.
About the assumptions to prove Point 3 of Thm.\ref{fix_point_3}, we note that existence of local sections for $\mG \ra \Omega$ is ensured under quite general conditions (Lemma \ref{lem_ls}); anyway, there are examples for which $\mG \ra \Omega$ does not admit local sections also in the case in which $X$ reduces to a single point: in explicit terms, there are compact groups $G$ with compact subgroups $H \subseteq G$ such that the projection $G \ra H \backslash G$ does not admit local sections (see \cite[I.7.5]{Ste}); but in this last case, we can use the Haar measure of $H$ and apply Lemma \ref{lem_hmh}, as done (implicitly) in \cite{DR86}.
Finally, unicity of $\mH_s$ is not ensured: by choosing a different section $s' : X \hra \Omega$ we may get a group bundle $\mH_{s'}$ not isomorphic to $\mH_s$; explicit examples are given in Sec.\ref{ss_cuntz}, by choosing $\Omega$ as in Sec.\ref{sec_nosec}.

%
%
%
%

Let $G$ be a compact Lie group. Then, the trivial bundle $\mG :=$ $X \times G$ admits an invariant $C(X)$-functional; moreover, for every group subbundle $\mH \subseteq \mG$ we find that $p_\Omega : \mG \ra \Omega$ has local sections (Lemma \ref{lem_ls}). Thus, we obtain the following result.

\begin{thm}
\label{cor_ds}
Let $G$ be a compact Lie group, $\mB$ a unital \sC algebra, and $\alpha : G \ra \aut \mB$ a strongly continuous action. Define $\mA := \mB^\alpha$, $\mZ := \mA' \cap \mA$, $C(\Omega) := \mB' \cap \mB$, $C(X) := C(\Omega) \cap \mA$ (so that, a bundle $\Omega \ra X$ is defined), and suppose that $\mA' \cap \mB = C(\Omega) \vee \mZ$. For every section $s \in$ $S_X(\Omega)$, the following properties are satisfied:
\begin{enumerate}
\item there is a compact group bundle $\mH_s \ra X$, $\mH_s \subseteq X \times G$, with a fibred $C(X)$-system $( \mF_s , \mH_s , \beta )$;
\item there is a $C(X)$-epimorphism $\eta : \mB \ra \mF_s$, which is injective on $\mA$, and such that $\eta (\mA) = \mF_s^\beta$, $\eta (\mA)' \cap \mF_s =$ $\eta (\mZ)$.
\end{enumerate}
\end{thm}

\subsection{Proof of Thm.\ref{fix_point_3}.}

Let $s \in S_X(\Omega)$. We denote by 
$\phi_s :$ $\mC \ra$ $C(X)$, $\phi_s (c) :=$ $c \circ s$, the $C(X)$-epimorphism associated with $s$. Moreover, we consider the group bundle ${\mH_s} \subseteq$ $\mG$ defined by
\begin{equation}
\label{def_h1}
{\mH_s} := 
\left\{ 
h \in \mG \ : \ \alpha^* ( h \ , \ s \circ p (h) ) \ = \ s \circ p (h) 
\right\} \ ,
\end{equation}
\noindent and denote by $H_x := \mH_s \cap G_x$, $x \in X$, the fibres of $\mH_s$. In the sequel, we will denote by $y,y', \ldots$ (resp. $h,h', \ldots$) generic elements of $\mG$ (resp. $\mH_s$). Note that if $c \in \ker \phi_s$, then by definition of $\mH_s$ we find $\alpha^x_h (c_x) \in$ $\ker \wa \phi_s$, $h \in \mH$, $x := p(h)$; in fact, $\alpha^x_h \circ c_x (s(x)) =$ $c_x \circ \alpha^* (h,s(x)) =$ $c_x (s(x)) =$ $0$.

Let us now suppose that there is a left invariant positive $C(X)$-functional $\varphi : C(\mG) \ra C(X)$. By Lemma \ref{stru_c}, Lemma \ref{iso_c}, there exists an isomorphism $\mC \simeq$ $C(\HsG)$ with $\Omega_x \simeq$ $H_x \backslash G_x$, $x \in X$. In the sequel, we will identify $\Omega$ with the homogeneous bundle $\HsG$.
This implies that we have a group bundle $p_\Omega : \mG \ra \Omega$, defined according to (\ref{def_ph}). The above considerations supply a proof of Point 1 of Thm.\ref{fix_point_3}.

We now consider the closed $C(X)$-ideal $\ker \phi_s \mB \subseteq \mB$ generated by $\left\{ 
c b , \right.$ $c \in \ker \phi_s ,$ $\left. b \in \mB \right\}$, and denote by 
\begin{equation}
\label{def_algf}
\eta : \mB \ra {\mF_s} := \left. \mB \right/ ( \ker \phi_s \mB)
\end{equation}
\noindent the associated $C(X)$-epimorphism. We denote by $\tau_x : {\mF_s} \ra {\mF_{s,x}}$, $x \in X$, the evaluation epimorphisms associated with $\mF_s$. Note that if we regard $\mB$ as a $C(\Omega)$-algebra, then $\eta$ is the restriction of $\mB$ over the closed subset $s(X) \subset$ $\Omega$.
By construction, $\ker \phi_s \mB$ is a non-degenerate $\ker \phi_s$-bimodule, thus every $b \in \ker \phi_s \mB$ admits a factorization $b = c b',$ $c \in \ker \phi_s$, $b' \in \ker \phi_s \mB$. By definition of $\mH_s$, we find $\alpha^x_h (c_x \pi_x(b) ) =$ $\alpha^x_h (c_x) \ \alpha^x_h \circ \pi_x (b) \in$ $\ker \wa \eta$, $x := p(h)$. This implies that $\mH_s$ satisfies (\ref{def_st}), thus we have a fibred action
\[
\beta : \mH_s \times_X \wa \mF_s \ra \wa \mF_s 
\ \ : \ \ 
\wa \eta \circ \alpha ( h,\sigma ) = \beta ( h , \wa \eta (\sigma) ) 
\ , \ \sigma \in \wa \mB \ .
\]

\begin{lem}
\label{lem_a_eta}
$\eta$ is injective on $\mA$. 
\end{lem}

\begin{proof}
In order to prove that $\eta$ is injective on $\mA$, it suffices to verify that each $\eta_x : \mB_x \ra \mF_{s,x}$, $\eta_x :=$ $\wa \eta |_{\mB_x}$, $x \in X$, is injective on $\pi_x(\mA) \subseteq$ $\mB_x$.
For this purpose, note that $\mB$ is a $C(\Omega)$-algebra, thus it is easy to verify that each $\mB_x$, $x \in X$, is a $C(\Omega_x)$-algebra; for every $\omega \in \Omega$, we consider $\mB_{x,\omega} :=$ $\mB_x / [ C_\omega (\Omega_x) \mB_x ]$, and denote by $\sigma_\omega \in \mB_{x,\omega}$ the image of $\sigma \in \mB_x$ w.r.t. the quotient epimorphism $\mB_x \ra \mB_{x,\omega}$. It is clear that 
$\left\| \sigma \right\| =$
$\sup_{\omega \in \Omega_x} \left\| \sigma_\omega \right\|$.
Moreover, since $\Omega_x$ is a homogeneous space w.r.t. the action $\alpha^*$, we find 
\begin{equation}
\label{eq_supo}
\left\| \sigma \right\|
=
\sup_{y \in G_x} \left\| \sigma_{\omega(y)} \right\|
\end{equation}
\noindent where $\omega (y) :=$ $\alpha^* ( y , s(x) )$, $x := p(y)$. Note that $\omega(e_x) =$ $s(x)$, where $e_x$ is the identity of $G_x$.
Now, by definition of $\eta$ we have $\eta_x (\sigma) =$ $\sigma_{s(x)}$; for every $y \in G_x$, $\sigma^0 \in \pi_x(\mA)$, by using the fact that $\sigma^0$ is $\alpha^x$-invariant, we find
$\sigma^0 + C_{\omega(y)}(\Omega_x) \mB_x =$
$\sigma^0 + \alpha^x_{y^{-1}} (C_{s(x)}(\Omega_x) \mB_x) =$
$\alpha^x_{y^{-1}} (\sigma^0 +  C_{s(x)}(\Omega_x)$.
Thus, we obtain
$\left\| \sigma^0_{\omega(y)} \right\| =$ 
$\inf_{\sigma , c} \left\| \sigma^0 + c \sigma \right\| =$
$\inf_{\sigma , c'} \left\| \sigma^0 + \alpha^x_{y^{-1}} (c' \sigma) \right\| =$
$\inf_{\sigma , c'} \left\| \alpha^x_{y^{-1}} (\sigma^0 +  c' \sigma ) \right\| =$
$\inf_{\sigma , c'} \left\| \sigma^0 +  c' \sigma  \right\| =$
$\left\| \sigma^0_{s(x)} \right\| =$
$\left\| \eta_x(\sigma^0) \right\|$,
where $c \in \mB_x$, $c \in C_{\omega(y)}(\Omega_x)$, $c' \in C_{s(x)}(\Omega_x)$. The above computation shows that if $\eta_x(\sigma^0) = 0$, then $\sigma^0_{\omega(y)} = 0$ for every $y \in G_x$, so that the lemma is proved.
\end{proof}

Let us now consider the $C(X)$-algebra tensor product $\gof :=$ $C(\mG) \otimes_X \mF_s$, and the associated \sC bundle $Q' : \wa \mF_s^\mG \ra X$ with fibres $\mF_{s,x}^\mG =$ $C(G_x) \otimes \mF_{s,x}$, $x \in X$. By (\ref{eq_pb}), we identify $\mF_s^\mG$ with $S_X(\mG,\wa \mF_s)$, and define the following $C(X)$-subalgebra of $\mF_s^\mG$:
\begin{equation}
\label{def_f_point}
\gbof 
:= 
\left\{ 
F \in \mF_s^\mG  \ : \ 
F (hy) = \beta ( h , F(y) )
\ , \
h \in \mH_s , y \in \mG
\right\} 
\end{equation}
\noindent (note that we assume $p(y) = p(h)$). The \sC bundle associated with $\gbof$ will be denoted by $\wa \mF_{s,\beta} \ra X$.
Let $\wa \alpha :$ $\mG \times_X \wa \mF_s^\mG \ra$ $\wa \mF_s^\mG$ denote the right translation $G$-action (\ref{def_ta}), so that $( \gof , \mG , \wa \alpha )$ is a fibred $C(X)$-system. If $F \in \gbof$ then
$\wa \alpha^x_y \circ F (hy') =$
$F ( hy'y ) =$ 
$\beta^x_h \circ F (y'y) =$ 
$\beta^x_h \circ \wa \alpha^x_y \circ F (y')$.
Thus, $\wa \mF_{s,\beta}$ is stable w.r.t. the action $\wa \alpha$, and a fibred $C(X)$-system 
\begin{equation}
\label{def_ib}
( \ \gbof \ , \ \mG \ , \ \wa \alpha \ )
\end{equation}
\noindent is defined. We call (\ref{def_ib}) the {\em induced} $C(X)$-{\em system} associated with $(\mB,\mG,\alpha)$. Some elementary remarks, analogous to the one in \cite[p.71]{DR86}, show that the fixed-point algebras ${\gbof}^{\wa \alpha}$ and $\mF_s^\beta$ coincide, in the sense that elements of ${\gbof}^{\wa \alpha}$ are exactly those of the type $\wa F (y) :=$
$\tau_{p(y)} (F)$, $F \in \mF_s^\beta$.

%
%

\begin{prop}
\label{prop_isob}
There is an equivariant $C(X)$-isomorphism
\[
T : ( \mB , \alpha ) \ra ( \gbof , \wa \alpha )
\ \ , \ \
Tb (y) := \wa \eta \circ \alpha ( y, \pi_x(b) ) 
\ , \
\]
\noindent where $b \in \mB$, $y \in \mG$, $x := p(y)$.
\end{prop}

\begin{proof}
By continuity of $\alpha$, $\wa \eta$, and since $b$ defines a section of $\wa \mB$, we find that $Tb$ belongs to $S_X ( \mG , \mF_s )$.  
In order to economize in notations, we define $\sigma :=$ $\pi_x(b)$, so that $Tb (y) =$ $\eta_x \circ \alpha^x_y (\sigma)$. We compute
$Tb (hy) = $
$\eta_x \circ \alpha^x_h \circ \alpha^x_y (\sigma) = $
$\beta^x_h \circ \eta_x \circ \alpha^x_y (\sigma) = $
$\beta^x_h \circ Tb (y)$,
thus $Tb \in$ $\gbof$. Let now $\wa T : \wa \mB \ra \wa \mF_{s,\beta}$ be the morphism induced by $T$. Then, we may regard $\wa T \sigma$ as a map from $G_x$ into $\mF_{s,x}$, and
$[ \wa \alpha ( y , \wa T \sigma )] \ (y') =$
$\wa T \sigma (y'y) =$
$\eta_x \circ \alpha^x_{y'} (\alpha^x_y(\sigma)) =$
$[\wa T \circ \alpha^x_y (\sigma)] \ (y')$.
So that, $\wa \alpha ( y , \wa T \sigma ) =$ $\wa T \circ \alpha ( y , \sigma)$, and $T$ is equivariant.
We now verify that $T$ is injective. For this purpose, we recall the reader to the $C(\Omega_x)$-algebra structure of $\mB_x$, and make use of the notation introduced in the proof of Lemma \ref{lem_a_eta}. Let $y \in \mG$, $\omega (y) :=$ $\alpha^* ( y , s(x) )$; then, for every $\sigma \in \mB_x$, we find
$\alpha^x_y ( \sigma + C_{\omega(y^{-1})} (\Omega_x) \mB_x  ) =$
$\alpha^x_y(\sigma) + C_{s(x)}(\Omega_x) \mB_x$.
By recalling the definition of the fibre epimorphism $\mB_x \ra \mB_{x,\omega}$, we conclude that $\alpha^x_y$ defines an isomorphism from $\mB_{x,\omega(y^{-1})}$ onto $\mB_{x,s(x)}$. This implies 
$\left\| \sigma_{\omega(y^{-1})} \right\| =$
$\left\| (\alpha^x_y(\sigma))_{s(x)} \right\|$.
On the other side, $(\alpha^x_y(\sigma))_{s(x)} =$ $\eta_x \circ \alpha^x_y (\sigma) =$ $\wa T \sigma (y)$. Thus, we conclude
$\left\| \wa T \sigma (y)  \right\| =$
$\left\| \sigma_{\omega(y^{-1})}  \right\|$,
$y \in G_x$.
\noindent Now, every $\omega' \in \Omega_x$ is of the type $\omega' = \omega (y)$ for some $y \in G_x$ . Thus (by (\ref{eq_supo})) $\wa T \sigma = 0$ if and only if $\sigma = 0$, and $\wa T$ is injective (i.e., $T$ is injective).
We now prove that $T$ is surjective.
Let $F \in \gbof$. Since every $\eta_x$ is surjective, for every $y \in \mG$ there is $b_x \in \mB$ such that $Tb_x (y) =$ $\wa \eta \circ \alpha ( y , \pi_x (b_x)) = $ $F(y)$, $x := p(y)$.
Let $\eps > 0$; for every $y \in \mG$, we consider the open set
\begin{equation}
\label{rel_ug}
U_y := \left\{ 
y' \in \mG \ : \
\left\|  F (y') - Tb_x (y')  \right\| < \eps 
\right\} \ .
\end{equation}
\noindent Now, for every $h \in \mH_s$ we have $\beta^x_h (F(y)) = F(hy)$, $\beta^x_h ( Tb_x(y') ) =$ $Tb_x (hy')$, and this implies
%
%
$U_y = U_{hy}$. By construction, $\left\{ U_y \right\}_y$ is an open cover for $\mG$; we pick a finite subcover $\left\{ U_{y_k} \right\}_k$, with the corresponding $b_k :=$ $b_{p(y_k)} \in$ $\mB$ satisfying (\ref{rel_ug}). Let us now consider the projection $p_\Omega :$ $\mG \ra$ $\Omega$; since $U_{y_k} = U_{hy_k}$, $h \in \mH_s$, there is a partition of unit $\left\{ \lambda_k \right\} \subset$ $\mC$ such that $\left\{ \lambda_k \circ p_\Omega  \right\} \subset$ $C(\mG)$ is subordinate to $\left\{ U_{y_k} \right\}$. We define
$b_\eps :=$ $\sum_k \lambda_k b_k \in$ $\mB$. 
\noindent If $y \in \mG$, then
$\left\|   F (y) - Tb_\eps (y)   \right\| \leq$
$\sum_k \lambda_k (p_\Omega (y)) \left\| F(y) - Tb_k (y)  \right\| \leq$
$\sum_k \lambda_k (y) \eps = \eps$.
This proves that the image of $T$ is dense in $\gbof$.
\end{proof}

Let $e_x \in G_x$ denote the identity. By definition of $\eta$ and $T$, we find
\begin{equation}
\label{eq_eta_T}
\tau_x \circ \eta (b) = Tb (e_x) 
\ \ , \ \
x \in X \ \ .
\end{equation}

\begin{cor}[Point 2 of Thm.\ref{fix_point_3}]
The map $T$ induces a $C(X)$-algebra isomorphism $\mA \simeq $ ${\mF_s}^\beta$. Moreover, ${\mF_s}^\beta = \eta (\mA)$.
\end{cor}

\begin{proof}
Since $T$ is covariant, we find that every $\wa \alpha$-invariant element of $\gbof$ is of the type $Ta$ for some $a \in \mA$. But $Ta$ is $\wa \alpha$-invariant if and only if $Ta (y'y) = Ta (y')$ for every $y' \in \mG$; thus, if $h \in \mH_s$ then
$Ta (h) =$
$Ta (he_x) =$
$\beta^x_h (Ta (e_x)) =$
$Ta (e_x h) =$
$Ta (e_x)$.
In other terms, by (\ref{eq_eta_T}),
\[
\beta^x_h \circ \tau_x ( \eta (a) ) =  \tau_x ( \eta (a) ) 
\ \ , \ \
x \in X \ ,
\]
\noindent so that we conclude $\eta (\mA) \subseteq \mF_s^\beta$. On the other side, it is clear that every $F \in \mF_s^\beta$ defines an $\wa \alpha$-invariant element of $\gbof$, say $\wa F$, by defining $\wa F (y) :=$ $\tau_{p(y)}(F)$, $y \in \mG$; in fact, $\wa F$ being fiberwise constant, we have $\wa F (hy) =$ $\wa F (y) =$ $\tau_x (F)$, $x :=$ $p(y) =$ $p(hy)$, so that $\beta^x_h \circ \wa F (y) =$ $\wa F (y) =$ $\wa F (hy)$. Thus, $\wa F = Ta$ for some $a \in \mA$, i.e. $F = \eta (a)$. Since $\eta (\mA)$ is isomorphic to $\mA$ (Lemma \ref{lem_a_eta}), the corollary is proved.
\end{proof}

Now, $C(X)$ may be regarded as a \sC subalgebra of $\eta(\mZ)$ (recall that $\mF_s$ is a $C(X)$-algebra), and $\eta (\mC) =$ $\phi (\mC) =$ $C(X) \subseteq$ $\eta (\mZ)$. Thus, we find
\begin{equation}
\label{eq_t2}
\eta ( \mA' \cap \mB ) = \eta ( \mC \vee \mZ ) = \eta (\mZ) \ \ .
\end{equation}

\begin{lem}[Point 3 of Thm.\ref{fix_point_3}]
\label{fix_point_2}
Suppose that $p_\Omega : \mG \ra \Omega$ has local sections.
If $\mA' \cap \mB = \mC \vee \mZ$, then $\eta (\mA)' \cap {\mF_s} = \eta (\mZ) \simeq \mZ$.
\end{lem}

\begin{proof}
It is clear that $\eta (\mZ) \subseteq \eta(\mA)' \cap {\mF_s}$. On the converse, let $F \in \eta(\mA)' \cap {\mF_s}$. In order to prove that $F \in \eta (\mZ)$, it suffices to verify that $F = \eta (b)$ for some $b \in \mA' \cap \mB$. In fact, in such a case (\ref{eq_t2}) implies that $F \in \eta (\mZ)$.
Let us now consider a finite open cover $\left\{ W_k \right\}$ of $\Omega$ satisfying the properties of Lemma \ref{lem_lm}. We consider a partition of unity $\left\{ \lambda_k \right\} \subset$ $C(\mG)$ subordinate to $\left\{ p_\Omega^{-1} (W_k) \right\} \subset$ $\mG$, so that $\lambda_k (hy) =$ $\lambda_k (y)$, $h \in \mH_s$, $y \in \mG$.
We consider the bundle morphisms (\ref{def_lm}), and define the map
\[
Tb (y) := \sum_k \lambda_k (y) \cdot \beta (  \ \delta_k (y) , \tau_{p(y)}(F) \ )
\ \ , \ \ 
y \in \mG \ .
\]
\noindent Since $F$ defines a continuous section of $\wa \mF_s$, and by continuity of $\beta$, $\delta_k$, we find $Tb \in$ $S_X(\mG, \wa \mF_s)$. Moreover, by applying (\ref{eq_lm}), we find $Tb (hy) =$ $\beta_h(Tb(y))$; thus, $Tb \in$ $\gbof$, i.e. $b \in \mB$. Since $F \in \eta (\mA)' \cap \mF_s$, we find $Tb \in$ $T (\mA)' \cap \gbof =$ $T ( \mA' \cap \mB )$. By applying (\ref{eq_eta_T}), we find $\tau_x \circ \eta (b) =$ $Tb (e_x) =$ $\tau_x (F)$ for every $x \in X$, and this implies $F = \eta (b)$.
\end{proof}

\begin{lem}[Point 3 of Thm.\ref{fix_point_3}]
\label{lem_hmh}
Suppose that there is a right-invariant positive $C(X)$-functional $\delta : C(\mH_s) \ra C(X)$. If $\mA' \cap \mB = \mC \vee \mZ$, then $\eta (\mA)' \cap {\mF_s} =$ $\eta (\mZ) \simeq \mZ$.
\end{lem}

\begin{proof}
As for the previous lemma, it suffices to prove that for every $F \in \eta (\mA)' \cap \mF_s$ there is $b \in$ $\mA' \cap \mB$ with $\eta (b) =$ $F$. To be more concise, we define ${\mF_s}' :=$ $\eta (\mA)' \cap \mF_s$; it is clear that ${\mF_s}'$ is stable w.r.t. $\beta$. Let us consider the map $\wa F \in$ $C(\mH_s) \otimes_X {\mF_s}'$, $\wt F (h) :=$ $\beta^x_h \circ \tau_x (F)$, $x := p(h)$. Since $\mH_s \subseteq$ $\mG$ is closed, by the Tietze extension theorem for continuous bundles (\cite[Chp.10]{Dix}) we find that there is $b \in$ $C(\mG)\otimes_X {\mF_s}'$ such that $b |_{\mH_s} =$ $\wt F$. In particular, $b (e_x) =$ $\tau_x (F)$, $x \in X$. Let us now consider the $\mH_s$-invariant mean $m : {\mF_s}' \ra ({\mF_s}')^\beta$ induced by $\delta$; then, for every $x \in X$, $F' \in {\mF_s}'$, we find $\tau_x \circ m(F') =$ $\int_{H_x} \beta^x_h \circ \tau_x (F') d \mu_x(h)$ (see (\ref{eq_int})). We define the map 
\[
\wt b (y) 
:= 
\int_{H_x} (\beta^x_h)^{-1} \circ b (hy) \ d \mu_x(h)
\ \in \wa \mF_s
\ \ , \ \ 
y \in \mG
\ ;
\]
since $b \in C(\mG)\otimes_X {\mF_s}'$, we find that $\wt b \in$ $C(\mG) \otimes_X {\mF_s}' \subseteq$ $C(\mG) \otimes_X \mF_s$; moreover, since $d\mu_x$ is right-invariant, we find $\wt b (h'y) =$ $\beta^x_{h'} \circ \int (\beta^x_{hh'})^{-1} \circ b (hh'y) d \mu_x(h) =$ $\beta^x_{h'} \circ \wt b (y)$. This implies $\wt b \in$ $\gbof$, thus $\wt b = Tb'$ for some $b' \in \mB$. Since $\wt b$ takes values in ${\mF_s}'$, we find $b' \in \mA' \cap \mB$; finally, $Tb'(e_x) =$ $\wt b (e_x) =$ $\int (\beta^x_h)^{-1} \circ \wt F (h) d \mu_x(h) =$ $\int \tau_x (F) d \mu_x(y) =$ $\tau_x(F)$, so that (by (\ref{eq_eta_T})) we find $\eta (b') = F$, and the lemma is proved.
\end{proof}

\subsection{A class of examples.}
\label{ss_cuntz}

Let $d \in \bN$, $d \geq 2$, $G \subseteq \sud$ a compact Lie group. We consider an Abelian $C(X)$-algebra $C(\Omega)$ carrying a strongly continuous action $\rho : G \ra \aut_X C(\Omega)$, such that $C(X) = C(\Omega)^\rho$. Starting from this data, we construct a \sC dynamical system $( \mB , G , \alpha )$ satisfying the hypothesis of Thm.\ref{cor_ds}.

We now start our construction. By Lemma \ref{stru_c}, every fibre $\Omega_x := q^{-1}(x)$ is homeomorphic to a homogeneous space $H_x / G$, $x \in X$. In order to simplify our notation, we write $\omega g :=$ $\omega \circ \rho_g$, $\omega \in \Omega$, $g \in G$.
Let us now consider the order $d$ Cuntz algebra $\mO_d$, generated by a set $\left\{ \psi_i \right\}_{i=1}^d$ of $d$ orthogonal isometries satisfying
\begin{equation}
\label{eq_cuntz}
\psi_i^* \psi_j = \delta_{ij} 1
\ \ , \ \ 
\sum_i \psi_i \psi_i^* = 1 
\ .
\end{equation}
\noindent  It is well-known that $G$ acts by automorphisms on $\mO_d$: if $g \in G$, then there exits a unique $\wa g \in {\bf aut} \mO_d$ such that
\begin{equation}
\label{def_g}
\wa g (\psi_i) = g \psi_i := \sum_{ij} g_{ij} \psi_j \ \ ,
\end{equation}
\noindent where $g_{ij} \in \bC$ are the matrix coefficients of $g$ w.r.t. the canonical basis of $\bC^d$. We define $\mB :=$ $C(\Omega) \otimes \mO_d$. $\mB$ is a trivial continuous bundle of Cuntz algebras over $\Omega$; moreover, $\mB$ is also a $C(X)$-algebra, with fibres $\mB_x \simeq$ $C (\Omega_x) \otimes \mO_d$, $x \in X$. In the sequel, we will regard elements of $\mB$ as continuous maps from $\Omega$ into $\mO_d$. We introduce the action
\[
\alpha :  G \ra {\bf aut}_X \mB
\ \ , \ \ 
[\alpha_g (b)]  (\omega) := \wa g \ ( b (\omega g^{-1}) )
\]
\noindent where $b \in \mB$, $g \in G$, $\omega \in \Omega$. We can now regard $\alpha$ as a fibred action $\alpha :$ $\mG \times_X \wa \mB \ra$ $\wa \mB$, with $\mG := X \times G$. We define $\mA := \mB^\alpha$; by construction
\[
\mA = \left\{ 
A \in \mB :  \wa g (A(\omega)) = A (\omega g)
\ \ , \ \ 
g \in \sud
\right\} \  .
\]

\begin{lem}
$\mA' \cap \mB = C(\Omega)$, and $\mA' \cap \mA = C(X)$.
\end{lem}

\begin{proof}
It is clear that $C(\Omega) \subseteq \mA' \cap \mB$.
Let $\mO_G$ denote the fixed-point algebra of $\mO_d$ w.r.t. the action (\ref{def_g}). By \cite[Cor.3.3]{DR87}, we conclude that $\mO'_G \cap \mO_d = \bC 1$. 
This also implies that $\mO'_G \cap \mB = C(\Omega)$ (we identify $\mO_G$ with the set of constant maps from $\Omega$ into $\mO_G$).
Let us now consider the constant map $A (\omega) := a \in \mO_G$, $A \in \mB$. Then
$[\alpha_g (A)] (\omega) =$
$ \wa g (a) = a = A (\omega)$, 
$g \in G$, $\omega \in \Omega$. This proves that $\mO_G \subseteq \mA$, so that $\mA' \cap \mB$ $\subseteq$ $\mO'_G \cap \mB$ $=$ $C(\Omega)$. We conclude that $\mA' \cap \mB = C(\Omega)$. In particular, if $b \in C(\Omega)$ is $\alpha$-invariant (i.e., $b \in \mA' \cap \mA$), then by definition $b$ is $\rho$-invariant. We conclude that $b \in C(\Omega)^\rho = C(X)$.
\end{proof}

{\em We now assume that there is a section} $s \in S_X(\Omega)$, and consider the associated $C(X)$-morphism $\phi_s :$ $C(\Omega) \ra$ $C(X)$. By applying Thm.\ref{cor_ds}, we obtain a fibred $C(X)$-system $( \mF , \mH_s , \beta )$, where $\mH_s \subseteq X \times G$ is defined according to (\ref{def_h1}). We now verify that for every $h \in \mH_s$, $F \in C(X) \otimes \mO_d$, $x := p (h)$, the following relations hold:
\begin{equation}
\label{eq_ex}
\left\{
\begin{array}{ll}
\mF \simeq C(X) \otimes \mO_d
\\
\beta^x_h \circ F (x) = \wa h \ (F(x)) 
\\
\mA = \mF^\beta 
\\
\mA \cap \mA' = \mA' \cap \mF  = C(X)
\end{array}
\right.
\end{equation}
\noindent About (\ref{eq_ex}.1): we have $\ker \phi_s :=$ $\left\{ c \in C(\Omega) : \right.$ $c (s(x)) = 0 ,$ $\left. x \in X \right\}$, thus a continuous map $b : \Omega \ra \mO_d$ belongs to $\ker \phi_s \mB$ if and only if $b \circ s(x) =$ $0$, $x \in X$. The quotient of $\mB$ w.r.t. such an ideal is the set of continuous maps from $s(X)$ into $\mO_d$, and $\eta$ is the restriction on $s(X)$. Since $s(X)$ is homeomorphic to $X$, we find (\ref{eq_ex}.1).
In order to simplify our notation, in the following lines we identify $\mF$ with $C(X) \otimes \mO_d$. 
About (\ref{eq_ex}.2): if $b \in \mB$, then $(\eta(b)) \ (x) =$ $b \circ s (x)$, $x \in X$. By definition of $\mH_s$, we find $s(x) h = h$ for every $h \in$ $H_x$, $x \in X$. Moreover, the fibred action $\beta$ satisfies $\beta^x_h \circ \eta_x = \eta_x \circ \alpha^x_h$, $h \in \mH$, $x := p(h)$. If $F \in \mF$ and $F = \eta (b)$ for some $b \in \mB$, then
$\beta^x_h \circ F (x)          =$
$\beta^x_h \circ b \circ s (x)  =$
$\alpha^x_h (b (s (x)) )        =$
$\wa h ( b (s(x)h^{-1}) )       =$
$\wa h ( b \circ s(x)   )       =$
$\wa h ( F(x) )$.
The equality (\ref{eq_ex}.3) follows from Thm.\ref{cor_ds}. 
About (\ref{eq_ex}.4), note that (\ref{eq_ex}.2) implies $C(X) \otimes \mO_G$ $\subseteq$ $\mF^\beta$. Thus, $(\mF^\beta)' \cap \mF$ $\subseteq$ $( C(X) \otimes \mO_G )'$ $\cap$ $\mF$; since $\mO'_G$ $\cap$ $\mO_d$ $=$ $\bC1$, we conclude that $( C(X) \otimes \mO_G )'$ $\cap$ $\mF$ $=$ $C(X)$, so that $(\mF^\beta)'$ $\cap$ $\mF$ $=$ $C(X)$. This also implies $(\mF^\beta)'$ $\cap$ $\mF^\beta$ $=$ $C(X)$.

Note that $\mH_s$ may be not full; in particular, the stabilizer $H$ of $\ker \eta$ in $G$ may be trivial (Rem.\ref{rem_fa_q}). Moreover, a different $s' \in S_X(\Omega)$ may define a group bundle $\mH_{s'} \subseteq \mG$ not isomorphic to $\mH_s$ (see Sec.\ref{sec_nosec})
Of course, the case in which $\Omega$ does not admit sections is possible (see Sec.\ref{sec_nosec}, where actually an action $\bT \ra$ $\aut_X C(\Omega)$ is considered; anyway, we may take $G := U(\bT)$, where $U : \bT \hra \sud$ is a faithful representation). In this case, the group bundle $\mH_s$ cannot be defined.

\section{Superselection structures.}
\label{sec_rc}

It is well-known that the field algebra describing a superselection structure of quantum localized observables is constructed in mathematical terms as a crossed product by a semigroup of \sC endomorphisms having permutation symmetry (\cite[\S 3]{DR89A},\cite[\S 2]{DR90}).
In the present section, we give a generalization of the above-cited crossed product to the nontrivial centre case, at least for the case of a single endomorphism (Thm.\ref{thm_dr}). In particular, we also cover the case of Hilbert \sC systems considered in \cite{BL03} (Thm.\ref{thm_bl}).

Let $\mA$ be a unital \sC algebra, $\rho \in {\bf end} \mA$ a unital endomorphism. We define $\mZ := \mA' \cap \mA$, and $C(X) := \left\{ f \in \mZ : \rho (f) = f \right\}$. Let us assume that there is a \sC algebra $\mB$ with identity $1$, carrying an inclusion $\mA \subset \mB$ of unital \sC algebras, such that:
\begin{enumerate}
\item  $\mB$ is generated as a \sC algebra by $\mA$ and a set $\left\{ \psi_i \right\}_{i=1}^d$, 
       $d \in \bN$, of isometries satisfying the Cuntz relations (\ref{eq_cuntz}). Thus, an 
       endomorphism $\sigma_\mB \in {\bf end} \mB$,
       \begin{equation}
       \label{def_sb}
       \sigma_\mB (b) := \sum_i \psi_i b \psi_i^* \ \ , \ \ b \in \mB
       \end{equation}
       \noindent is defined. By universality of the Cuntz algebra, there is a unital monomorphism 
       $j : \mO_d \hra \mB$;
\item  the following relations hold:
       \begin{equation}
       \label{eq_rho}
       \rho (a) = \sigma_\mB (a)
       \ \ , \ \
       a \in \mA \ .
       \end{equation}
       \noindent Note that the previous relations imply that if $f \in C(X)$, then $f$ commutes with 
       $\psi_1$, $\ldots$, $\psi_d$, thus $f \in$ $\mB' \cap \mB$. In other words, $\mB$ is a 
       $C(X)$-algebra (also note that $\sigma_\mB$ is a $C(X)$-endomorphism);
\item  there is a strongly continuous action $\alpha : \sud \ra {\bf aut}_X \mB$, such that $\mB^\alpha = \mA$, $\alpha_g \circ j (t) = j \circ \wa g (t)$, $g \in \sud$, $t \in \mO_d$, where $\wa g \in {\bf aut} \mO_d$ is defined as in (\ref{def_g}). 
\end{enumerate}

We denote by $\mC$ the centre of $\mB$, and by $\Omega$ the spectrum of $\mC$.

Let $\mO_\sud \subset \mO_d$ be the fixed-point algebra w.r.t. the $\sud$-action (\ref{def_g}); then, the above considerations imply that $j(\mO_\sud) \subseteq \mA$. For every $r,s \in \bN$, we consider the intertwiner spaces 
\[
\left\{
\begin{array}{ll}
( \rho^r , \rho^s ) := 
\left\{  
t \in \mA : \rho^s(a)t = t \rho^r(a) , a \in \mA 
\right\}
\\
( \sigma_\mB^r, \sigma_\mB^s ) := 
\left\{  
w \in \mB : \sigma_\mB^s(b)w = w \sigma_\mB^r(b) , b \in \mB 
\right\}
\end{array}
\right.
\]
\noindent In particular, for $r = 0$, we define $\sigma_\mB^0 := \iota_\mB$ (the identity on $\mB$) and $\rho^0 := \iota_\mA$ (the identity on $\mA$). Let $I := \left\{ i_1 , \dots , i_r  \right\}$ be a multiindex of lenght $|I| = r \in \bN$. We define $\psi_I := \prod_r \psi_{i_r} \ \in ( \iota_\mB , \sigma_\mB^r )$; note that 
\begin{equation}
\label{eq_cI}
\sum_I \psi_I \psi_I^* = 1 \ \ , \ \ \psi_I^* \psi_{I'} = \delta_{II'} 1 \ . 
\end{equation}
\noindent Let $J$ be a multiindex with lenght $s \in \bN$. By \cite[\S 2]{DR87}, the following "flip" operators belong to $j(\mO_\sud) \subseteq \mA$:
\[
\eps (r,s) 
\ := \
\sum_{IJ} \psi_J \psi_I \psi_J^* \psi_I^* 
\ \in \ 
( \sigma_\mB^{r+s} , \sigma_\mB^{r+s} ) \ ;
\]
\noindent note that in particular $\eps (r,s) \in ( \rho^{r+s} , \rho^{r+s} )$. To be coherent, we also define $\eps (0,s) = \eps (r,0) = 1$. By using (\ref{eq_cI}), we compute
\begin{equation}
\label{eq_sIe}
\sigma_\mB (\psi_I) =
\sum_i \psi_i \psi_I \psi_i^* =
\sum_{iJ} \psi_i \psi_J ( \psi_J^* \psi_I ) \psi_i^* =
\eps (1,r) \psi_I \ \ .
\end{equation}
\noindent With the same argument used for the $\eps(r,s)$'s, we conclude that if
 $\bP(r)$ is the permutation group of $r$ objects and 
\[
\eps (p) := \sum  \psi_{i_1} \cdots \psi_{i_r}  
                  \psi_{i_{p(r)}}^* \cdots \psi_{i{p(1)}}^*  \ ,
\]
\noindent then $\eps(p) \in ( \rho^r , \rho^r  )$. 
Let us now define
\[
( \rho^r , \rho^s  )_\eps
:= 
\left\{ \ 
t \in ( \rho^r , \rho^s  )
:
\rho (t) \eps (r,1) = \eps(s,1) t
\ \right\}
\ .
\]
\noindent By definition of $\eps (r,s)$ ($r=s=0$), it turns out $( \iota_\mA , \iota_\mA )_\eps = C(X)$; moreover, a simple computation shows that $\eps (r,s)$ $\in$ $( \rho^{r+s} , \rho^{r+s} )_\eps$, $r,s \in \bN$. Note that every $(\rho^r , \rho^s)_\eps$ is a $C(X)$-bimodule w.r.t. left and right multiplication by elements of $C(X)$. We say that $\rho$ has {\em permutation quasi-symmetry} if every $t \in ( \rho^r , \rho^s )$ admits a decomposition $t = \sum_i \rho^s (z_i) t_i$, $z_i \in \mZ$, $t_i \in ( \rho^r , \rho^s )_\eps$ (note that since $t_i \in ( \rho^r , \rho^s  )$, then $\rho^s(z_i)t_i = t_i \rho^r (z_i)$). In particular, if $( \rho^r , \rho^s ) = ( \rho^r , \rho^s )_\eps$, $r,s \in \bN$, then we say that $\rho$ has {\em permutation symmetry}, according to \cite[\S 4]{DR89A}.

\begin{prop}
\label{prop_rc}
If $\rho$ has permutation quasi-symmetry, then $\mA' \cap \mB = \mC \vee \mZ$.
\end{prop}

\begin{proof}
We proceed in a order of ideas similar to \cite[Lemma 5.1]{DR88}. As a first step, note that $\mA' \cap \mB$ is $\sud$-stable; thus, by Fourier analysis, $\mA' \cap \mB$ is generated as a \sC algebra by sets of "irreducible tensors" of the type $\left\{ T_1 , \dots , T_n \in \mA' \cap \mB \right\}$, satisfying the relations 
\begin{equation}
\label{irr_tensors}
\alpha_g (T_i) = \sum_j^N T_j \ u_{ji} (g) \ , 
\end{equation}
\noindent where $g \in \sud$ and $u_{ji} (g) \in \bC$ are matrix element in some irreducible, unitary representation of $\sud$ (in particular, note that $\ovl{u_{ji} (g)} = u^*_{ij} (g)$ are the matrix elements of the inverse of $(u_{ij}(g))_{ij}$). Now, every irreducible representation of $\sud$ is a subrepresentation of some tensor power of the defining representation. It follows from (\ref{def_g}) that $\sud$ acts on ${\mathbb{H}}_d := $ ${\mathrm{span}} \left\{ \psi_i  \right\}_i$ $\subset$ $( \iota_\mB , \sigma_\mB  )$ as the defining representation. Moreover, for every $r \in \bN$ we can identify ${\mathbb{H}}_d^r$ $:=$ ${\mathrm{span}} \left\{ \psi_I \in  ( \iota_\mB , \sigma_\mB^r )  , |I| = r  \right\}$ with the $r$-fold tensor power of ${\mathbb{H}}_d$, and $\sud$ acts on ${\mathbb{H}}_d^r$ as the $r$-fold tensor power of the defining representation. Thus, there is $n \in \bN$ and an orthonormal set $\left\{ \varphi_h \right\} \subset {\mathbb{H}}_d^n$ (i.e., $\varphi_i^* \varphi_j = \delta_{ij}1$), such that 
\begin{equation}
\label{irr_psi}
\alpha_g ( \varphi_i ) = \sum_j \varphi_j u_{ji} (g) \ .
\end{equation}
\noindent Let $W := \sum_i T_i \varphi_i^*$; by using (\ref{irr_tensors}) and (\ref{irr_psi}), we find
$\alpha_g (W) =$
$\sum_{ihk} T_h \varphi_k^* \ u_{hi}(g) \ \ovl{u_{ki}(g)} =$
$\sum_{hk} \delta_{hk} T_h \varphi_k^* =$
$W$.
We conclude that $W \in \mA$. Moreover, $T_i$ $=$ $W \varphi_i$, where $W \in ( \rho^n , \iota_\mA )$, $\varphi \subset (\iota_\mB , \sigma_\mB^n)$. Thus, by permutation quasi-symmetry, $W = \sum_k \rho^{pd}(z_k) W_k$, where $z_k \in \mZ$ and $W_k \in (\rho^n , \iota_\mA)_\eps$. So that, we obtain $T_i$ $=$ $\sum_k z_k W_k \varphi_i$. Let us now define $\ovl T_{i,k} := W_k \varphi_i$; we prove that $\sigma_\mB (\ovl T_{i,k}) = \ovl T_{i,k}$ (so that, $\ovl T_{i,k}$ commutes with $\psi_1$, $\ldots$, $\psi_d$ and $\ovl T_{i,k} \in \mC$). By permutation quasi-symmetry and (\ref{eq_sIe}), we obtain 
$\sigma_\mB (\ovl T_{i,k}) =$
$\rho (\ovl W_k) \sigma_\mB (\varphi_i) =$
$\ovl W_k \eps(1,n) \cdot \eps(n,1) \varphi_i =$
$\ovl T_{i,k}$,
and this implies $\ovl T_{i,k} \in \mC$. Since $T_i = \sum_n z_k \ovl T_{i,k}$, we conclude that $\mA' \cap \mB = \mC \vee \mZ$.
\end{proof}

We now look for \sC epimorphisms $\eta : \mB \ra \mF$ which are injective on $\mA$, and such that $\eta (\mA)' \cap \mF =$ $\eta (\mZ)$. Pairs $( \mF , \eta )$ of the above type are called {\em Hilbert extensions of} $\mA$. {\em Every Hilbert extension corresponds to the crossed product of $\mA$ by $\rho$ in the sense of \cite[\S 4]{DR89A}, thus can be interpreted as a "field algebra" associated with} $(\mA,\rho)$. We discuss existence and unicity of the Hilbert extension in two important cases.

\

{\bf Doplicher-Roberts endomorphisms (with non-trivial centre).}
Let $\mA$ be a unital \sC algebra with centre $\mZ$, and $\rho \in {\bf end} \mA$ an endomorphism satisfying the special conjugate property in the sense of \cite[\S 4]{DR89A} (i.e., $\rho$ has permutation symmetry, and there is an isometry $R \in (\iota_\mA,\rho^d)_\eps$ satisfying the {\em special conjugate equations} $R^* \rho (R) =$ $(-1)^{d-1} d^{-1} 1$, $RR^* =$ ${d!}^{-1} \sum_{p \in \bP_d} \mathrm{sign}(p) \eps (p)$). 
Then, we can construct a \sC algebra $\mB$ satisfying the above properties, as the crossed product of $\mA$ by the dual $\sud$-action induced by $\rho$ (\cite[Thm.4.2]{DR89A}); we maintain the notation $\mC \equiv C(\Omega)$ for the centre of $\mB$, and $C(X)$ for the \sC algebra of $\rho$-invariant elements of $\mZ$. Since in particular $\rho$ has permutation symmetry the previous proposition applies, so that $\mA' \cap \mB =$ $\mC \vee \mZ$.

\begin{thm}
\label{thm_dr}
Let $\rho \in {\bf end} \mA$ be an endomorphism satisfying the special conjugate property. Sections $s \in$ $S_X(\Omega)$ are in one-to-one correspondence with Hilbert extensions $(\mF_s,\eta_s)$ of $\mA$. For every $s \in S_X(\Omega)$, there exists a group bundle $\mG_s \subseteq$ $X \times \sud$ such that $( \mF_s , \mG_s , \beta )$ is a fibred $C(X)$-system satisfying $\eta_s (\mA) = {\mF_s}^\beta$, $\eta_s (\mA)' \cap \mF_s =$ $\eta_s (\mZ)$.
\end{thm}

\begin{proof}
If there is $s \in S_X(\Omega)$, then existence of a Hilbert extension $( \mF_s , \eta_s )$ follows from Thm.\ref{cor_ds}, and a fibred action $( \mF_s , \mG_s , \beta )$ is defined, satisfying the required properties. On the other side, if $( \mF , \eta )$ is a Hilbert extension, then we note that $\eta$ is injective on $j(\mO_d) \subset \mB$, in fact $\mO_d$ is simple and $\eta (t) \neq 0$ for every $t \in j(\mO_d) \cap \mA$. Let us now consider the restriction $\phi := \eta |_\mC$. Then, $\phi (\mC) \subseteq \eta (\mZ)$, thus for every $c \in \mC$ there is $f \in \mZ$ with $\phi (c) =$ $\eta (f)$. Since $\eta$ is injective on $\mZ$, we find that $f$ is unique. Moreover, by using the fact that $\phi (c) \in \mF' \cap \mF$, we find $\eta \circ \rho (f) =$ $\eta \circ \sigma_\mB (f) =$ $\sum_i \eta (\psi_i) \eta (f) \eta (\psi_i)^* =$ $\eta(f)$. We conclude that $\eta ( f - \rho (f) ) = 0$, i.e. $f = \rho(f)$; by definition of $C(X)$, this implies $f \in C(X)$. Thus, we have a \sC epimorphism $\phi : \mC \ra \eta (C(X))$, which is injective on $C(X) \subseteq$ $\mC$. By applying the Gel'fand transform, we obtain the desired section $s \in$ $S_X(\Omega)$.
\end{proof}

Note that $\mG_s$ may be a non-trivial bundle, and may depend on the choice of $s$ (see Sec.\ref{sec_nosec}). The question of unicity of $\mG_s$ may be approached by regarding at $\mF_s$ as a crossed product $\mA \rtimes \wa \mG_s$ in the sense of \cite{Vas05}. Anyway there are some further complications, due to the fact that we may choose a Hilbert $C(X)$-bimodule instead of ${\mathbb{H}}_d :=$ ${\mathrm{span}}$ $\left\{ \psi_i \right\}$ to construct an analogue of the \sC algebra $\mB$ (see \cite[\S 1]{Vas03}). For this reason, we postpone a complete discussion to a forthcoming paper.

%
%

\

{\bf Baumg\"artel-Lled\'o endomorphisms.}
The following class of \sC endomorphisms has been studied in \cite[\S 4]{BL03}.
Let $\mA$ be a unital \sC algebra;
we say that $\rho$ is a {\em canonical endomorphism} of $\mA$ if $\rho$ has permutation quasi-symmetry, with the additional property that every $(\rho^r , \rho^s)_\eps$ is a free $C(X)$-bimodule generated by a finite-dimensional vector space, say $( \rho^r , \rho^s )_\bC$. The family of vector spaces $( \rho^r , \rho^s  )_\bC$, $r,s \in \bN$, is required to be $\rho$-stable, closed for multiplication (i.e., $( \rho^k , \rho^s  )_\bC ( \rho^r , \rho^k  )_\bC \subseteq$ $( \rho^r , \rho^s  )_\bC$) and involution, and such that $(\iota_\mA , \iota_\mA)_\bC =$ $\bC$, $\eps (r,s) \in (\rho^{r+s} , \rho^{r+s})_\bC$.  
Let us now suppose that there is $R \in$ $( \iota , \rho^d )_\eps$ satisfying the special conjugate equations; in such a case, we say that $\rho$ is a {\em special canonical endomorphism}. The $\rho$-stable \sC subalgebra of $\mA$ generated by $\left\{ R , \eps (r,s) \right\}_{r,s}$, is isomorphic to $\mO_\sud$ (\cite[Thm.4.1]{DR87}): such an isomorphism defines a dual action $\mu : \mO_\sud \hra \mA$. Again, $\mB$ can be constructed as the crossed product of $\mA$ by $\mu$ in the sense of \cite[\S 3]{DR89A}, and Prop.\ref{prop_rc} applies.

\begin{lem}
There is a closed group $G \subseteq \sud$ with a bundle isomorphism $\Omega \simeq X \times G \backslash \sud$. The group $G$ is unique up to conjugation in $\sud$.
\end{lem}

\begin{proof}
We retain the same notation used in the proof of Prop.\ref{prop_rc}. As a preliminary step, we consider the \sC subalgebra $\mC_0$ of $\mB$ generated by the set $\left\{ W^* \varphi \right. :$ $W \in ( \iota_\mA , \rho^n  )_\bC ,$ $\left. \varphi \in {\mathbb{H}}_d^n  \right\}$. Since ${\mathbb{H}}_d^n \subset ( \iota_\mB , \sigma_\mB^n  )$, we find that $W^* \varphi a = W^* \rho^n(a) \varphi = a W^* \varphi$, $a \in \mA$. Moreover, by using (\ref{eq_sIe}), we find $\sigma_\mB (W^* \varphi) = W^* \eps (1,n) \eps (n,1) \varphi = W^* \varphi$: this means that $W^* \varphi$ commutes with $\psi_i$ for every $i = 1 , \ldots , d$ (recall (\ref{def_sb})). We conclude that $\mC_0$ is contained in $\mC$. Morever, $W_1^* \varphi_1 W_2^* \varphi_2 =$ $W_1^* \rho^{n_1} (W_2) \varphi_2 \varphi_1 =$ $(W_1 W_2)^* \varphi_1 \varphi_2$, for every $W_1^* \varphi_1$, $W_2^* \varphi_2$ $\in \mC_0$: thus, the set of linear combinations of terms of the type $W^* \varphi$ (and the adjoints $\varphi^* W$) is dense in $\mC_0$. Let now $c := \sum_i W_i^* \varphi_i$ such that $\alpha_g (c) = c$, $g \in \sud$; by averaging w.r.t. the Haar mean, say $m : \mB \ra \mA$, we find $c = m (c) = \sum_i W_i^* m(\varphi_i)$. Now, $m(\varphi_i)$ is an $\alpha$-invariant element of some ${\mathbb{H}}_d^{n_i}$. Since $\alpha_g \circ j = j \circ \wa g$, $g \in \sud$, we conclude that $m(\varphi_i)$ belongs to $j (\mO_\sud) \subseteq \mA$. Now, $j (\mO_\sud)$ is generated by elements belonging to the family $( \rho^r , \rho^s  )_\bC$, $r,s \in \bN$, thus $m(\varphi_i) \in ( \iota_\mA , \rho^{n_i} )_\bC$. This implies that $c =$ $m(c) =$ $\sum_i W_i^* m(\varphi_i)$ belongs to $\bC 1$, in fact $W_i^* m(\varphi_i) \in$ $( \rho^{n_i} , \iota_\mA  )_\bC \cdot ( \iota_\mA , \rho^{n_i} )_\bC \subseteq ( \iota_\mA , \iota_\mA  )_\bC = \bC 1$. This proves that $\mC_0^\alpha \simeq \bC 1$; thus, by well-known results, if we pick an element $\omega_0$ of the spetrum $\Omega_0$ of $\mC_0$, and define $G \subseteq \sud$ as the stabilizer of $\omega_0$, then we get a homeomorphism $\Omega_0 \simeq G \backslash \sud$. Note that $G$ is unique up to conjugation in $\sud$.
We now introduce the map $\Gamma : C(X) \otimes \mC_0 \ra \mC$, $\Gamma ( f \otimes c_0 ) := fc_0$. Let $W^* \varphi \in \mC_0$; then, $\left\| f W^* \varphi  \right\|^2 = \left\| (f^*f) ( \varphi^* WW^* \varphi  )  \right\|$. Now, note that $c_{\varphi,W}:= ( \varphi^* WW^* \varphi )$ belongs to $( \iota_\mA , \iota_\mA )_\bC$; thus, $c_{\varphi,W}$ is a (positive) multiple of the identity, and it is now clear that $c_{\varphi,W} = \left\| W^* \varphi  \right\|^2$. Moreover, $\left\| f W^* \varphi  \right\|^2 =$ $\left\| f^*f c_{\varphi,W}  \right\|^2 =$ $\left\| f  \right\|^2 c_{\varphi,W} =$ $\left\| f  \right\|^2 \left\| W^* \varphi  \right\|^2$. This proves that $\Gamma$ is isometric. We now prove that $\Gamma$ is surjective. With the same argument as the previous proposition, $\mC$ is the closed vector space spanned by multiplets of the tipe $\left\{ c_i \right\}_{i=1}^n$, such that $\alpha_g (c_i) = \sum_k c_k u_{ki} (g)$, $g \in \sud$; the coefficients $u_{ki} (g) \in \bC$ are matrix elements of some irreducible representation of $\sud$, say $u$. Thus, there is a set $\left\{ \varphi_i  \right\} \subset ( \iota , \sigma_\mB^n )$, $\varphi_i^* \varphi_j = \delta_{ij} 1$ such that $W^* := \sum_i c_i \varphi_i^*$ belongs to $\mA$ (i.e., $\alpha_g(W^*) = W^*$); note that $c_i = W^* \varphi_i$. Moreover, it is clear that $W^* \rho^n(a) = a W^*$, i.e. $W \in (  \iota_\mA , \rho^n  )$. Since $\rho$ is a canonical endomorphism, we find $W = \sum_h f_h W_h$, $W_h \in ( \iota , \rho^n  )_\bC$. Thus, $c_i = \sum_{kh} f_h^* ( W_h^* \varphi_i )$, and $\Gamma$ is surjective. It is now clear that $\Omega \simeq$ $X \times \Omega_0 \simeq$ $X \times G \backslash \sud$.
\end{proof}

\begin{thm}
\label{thm_bl}
Let $\rho \in {\bf end} \mA$ be a special canonical endomorphism. Then, there exists a Hilbert extension $( \mF , \eta )$, and a closed group $G \subseteq \sud$ unique up to conjugation in $\sud$, with an action $\beta : G \ra \aut \mF$ such that $\mA \simeq$$\eta (\mA) = \mF^\beta$, $\eta (\mA)' \cap \mF =$ $\eta (\mZ)$. 
\end{thm}

\begin{proof}
It suffices to apply Thm.\ref{cor_ds} to the \sC dynamical system $( \mB , \sud )$, by choosing $s \in S_X(X \times (G \backslash \sud))$ of the type $s (x) := ( x , \omega_0 )$, $\omega_0 \in G \backslash \sud$.
\end{proof}

%
%

\section{Appendix.}
\label{sec_nosec}

{\bf Bundles of homogeneous spaces lacking of sections.}
Let us denote by $\bT := \left\{ z \in \bC : |z| = 1 \right\}$ the torus. For every $n \in \bN$, we consider the compact subgroup 
$R_n := \left\{ z \in \bT : z^n = 1  \right\}$ 
of roots of unity. Elementary computations show that the coset space $R_n \backslash \bT$ 
can be identified with $\bT$, and that there is an exact sequence 
${\bf 1} \hra$
$R_n \stackrel{\subset}{\lra}$
$\bT \stackrel{p_n}{\lra}$
$\bT \ra$
${\bf 1}$,
where ${\bf 1}$ is the trivial group, and $p_n (z) := z^n$, $z \in \bN$.

Let $S^2$ denote the $2$-sphere. Principal $\bT$-bundles over $S^2$ are classified by the cohomology group $H^1(S^2,\bT) \simeq \bZ$: for every $h \in \bZ$, we denote by $\Omega_h \ra S^2$ the associated principal $\bT$-bundle, having fibre $\Omega_{h,x} \simeq \bT$. The bundle $\Omega_h \ra S^2$ is trivial if and only if $h = 0$. Let $k = 1 , \ldots , n$; it follows from \cite[Prop.7.1.7]{Hus} that there is an action
\[
\rho^k : \bT  \ra  {\bf aut}_X C(\Omega_h)
\ \ , \ \
\]
\noindent which fiberwise behaves as the multiplication $\rho^{k,x} :$ $\bT \ra$ $\aut C(\Omega_{h,x})$, $\rho^{k,x}_z c (\omega)$ $:=$ $c(\omega z^k)$, $c \in C(\Omega_{h,x})$, $\omega \in$ $\bT \simeq$ $\Omega_{h,x}$. Since multiplication by $z^k$, $z \in \bT$, defines an ergodic action on $\bT$, we conclude that the fixed-point algebra of $C(\Omega_{h,x})$ w.r.t. $\rho^{k,x}$ reduces to the complex numbers. Thus, $C(\Omega)^\rho \simeq C(S^2)$.
By general properties of principal bundles (\cite[Thm.6.2.3]{Hus}), the unique principal $\bT$-bundle over $S^2$ which admits a section is the trivial bundle $\Omega_0 \simeq$ $S^2 \times \bT$. 

%
%

\

{\bf Non-isomorphic group subbundles, with isomorphic associated homogeneous bundles.}
Let $n \in \bN$, and $\bO(n) \subset \bM_n(\bR)$ denote the orthonomormal group with unit $1_n$. For every $m < n$, 
we have that $\bO(m)$ may be regarded as a closed subgroup of $\bO(n)$ via the embedding $\bO(m) \oplus 1_{n-m} \subset$ $\bO(n)$. The quotient $\bO(m) \backslash \bO(n)$ is the homogeneous space known as the Stiefel manifold (\cite[7.1]{Hus}).

Let $X$ be a compact Hausdorff space. For every $n \in \bN$, we denote by $T_n :=$ $X \times \bR^n$ the trivial rank $n$ vector bundle, and by $\mG_n :=$ $X \times \bO (n)$ the trivial group bundle with fibre $\bO(n)$.
For every rank $d$ real vector bundle $\mE \ra X$, we consider the associated bundle of orthonormal endomorphsms $\mathcal{OE} \ra X$ with fibre $\bO (d)$ (\cite[I.4.8]{Kar}). It is well-known that if $\mE' \ra X$ is a real rank $d$ vector bundle, then $\mathcal{OE}$ is isomorphic to $\mathcal{OE}'$ if and only if $\mE' =$ $\mE \otimes \mL$ for some real line bundle $\mL \ra X$.

Thus, if $X$ is a space such that the cohomology group $H^1(X,\bZ_2)$ is trivial, then every real line bundle $\mL \ra X$ is trivial, and $\mE \simeq \mE'$ if and only if $\mathcal{OE} \simeq \mathcal{OE}'$. 
Let now $\mE \ra X$ be a non-trivial rank $d$ vector bundle such that $\mE \oplus T_n \simeq$ $T_{d+n}$ for some $n \in \bN$ (i.e., $\mE$ has trivial class in $K$-theory). Since $\mE$ is non-trivial, $\mathcal{OE}$ is non-trivial; on the other side, $\mathcal O T_{d+n} =$ $\mG_{d+n}$ is trivial. Moreover, the natural monomorphism $j : \mE \hra$ $T_{d+n}$ induces a monomorphism ${\bf ad} j : \mathcal{OE} \hra$ $\mG_{d+n}$. So that, the trivial bundle $\mG_{d+n}$ has two non-isomorphic group subbundles with fibre $\bO (d)$, i.e. the trivial bundle $\mG_d$ and $\mH :=$ ${\bf ad} j (\mathcal{OE})$. This happens despite the fact that we have isomorphisms
\begin{equation}
\label{eq_om}
\mG_d \backslash \mG_{d+n} 
\ \ \simeq \ \
\mH \backslash \mG_{d+n} 
\ \ \simeq \ \
\Omega := X \times ( \bO(d) \backslash \bO(d+n) )
\ ;
\end{equation}
\noindent in fact, $\Omega$ is a bundle having the (trivial) cocycle associated with $\mG_{d+n}$ as a set of transition maps (see \cite[Thm.6.4.1]{Hus} and following remarks). Now, $C(\Omega)$ is endowed with an automorphic action $\rho : \bO(d+n) \ra {\bf aut}_X C(\Omega)$ such that $C(\Omega)^\rho =$ $C(X)$, and $\Omega$ is clearly full (in fact, it is a trivial bundle over $X$). This implies that $\Omega \ra X$ is a homogeneous bundle. Anyway, (\ref{eq_om}) implies that $\Omega$ can be recovered as a quotient of $\mG_{d+n}$ w.r.t. non-isomorphic subbundles, namely $\mG_d$ and $\mH$.

For example, all the above considerations apply in the case in which $X$ is the sphere $S^k$, $k \neq 1,3,7$, and $\mE$ is the the tangent bundle $TS^k$ (\cite[I.5.5]{Kar}).

\

{\bf Acknowledgments.} The author would like to thank S. Doplicher, C. Pinzari for stimulating discussions, and Gerardo Morsella, for precious help.



{\small
\begin{thebibliography}{99}

%
%

%
%

\bibitem{Bla96}
E. Blanchard: D\'eformations de \sC alg\'ebres de Hopf, Bull. Soc. math. France {\bf 124}, 141-215 (1996).

%
%
%
%
%
%


\bibitem{BL03}
H. Baumg\"artel, F. Lled\'o: Duality of compact groups and Hilbert \sC systems for \sC algebras with a nontrivial center, Int. J. Math. 15, 759-812 (2004).


%
%


\bibitem{CG04}
A. Carey, H. Grundling: Amenability of the Gauge Group, Lett.Math.Phys. {\bf 68}, 113-120 (2004).


\bibitem{Cio03}
F. Ciolli: Massless scalar free Field in 1+1 dimensions I: Weyl algebras Products and Superselection Sectors, preprint arXiv math-ph/0511064 (2005).

%
%

\bibitem{Dix}
J. Dixmier: \sC {\it algebras}, North-Holland Publishing Company, Amsterdam - New
York . Oxford (1977).

%
%
%
%
%
%
%


\bibitem{DR86}
S. Doplicher, J. E. Roberts: A Remark on Compact Automorphism Groups of \sC Algebras, J. Funct. Anal. {\bf 66}, 67-72 (1986).


\bibitem{DR87}
S. Doplicher, J.E. Roberts: Duals of Compact Lie Groups Realized in the Cuntz Algebras and Their Actions on \sC Algebras, J. Funct. Anal. {\bf 74} 96-120 (1987).


\bibitem{DR88}
S. Doplicher, J.E. Roberts: Compact Group Actions on \sC Algebras, J. Operator Theory {\bf 91} 227-284 (1988).


%
%


\bibitem{DR89A}
S. Doplicher, J.E. Roberts: Endomorphisms of \sC algebras, Cross Products and Duality for Compact Groups, Annals of Mathematics {\bf 130}, 75-119 (1989).


\bibitem{DR90}
S. Doplicher, J.E. Roberts: Why there is a field algebra with a compact gauge group describing the superselection structure in particle physics, Commun. Math. Phys. {\bf 131} 51-107 (1990).



%
%
%
%

%
%

\bibitem{Gie82}
G. Gierz: Bundles of topological vector spaces and their duality, Lectur Notes in mathematics, 955, Springer-Verlag (1982).



\bibitem{Hir}
F. Hirzebruch, {\it Topological Methods in Algebraic Geometry}, Springer-Verlag, 1966.

\bibitem{Hus}
D. Husemoller: {\it Fiber Bundles}, Mc Graw-Hill Series in Mathematics, 1966.

\bibitem{Kar}
M. Karoubi: {\it K-Theory}, Springer Verlag, Berlin - Heidelberg - New York, 1978.

\bibitem{Kas88}
G.G. Kasparov, Equivariant $KK$-Theory and the Novikov Conjecture, Invent. Math. {\bf 91}, 147-201 (1988).

%
%

\bibitem{MS90}
G. Mack, V. Schomerus: Conformal field algebras with quantum symmetry from the theory of superselection sectors, Commun. Math. Phys. {\bf 134} 139-196 (1990).


\bibitem{Nil96}
M. Nilsen: \sC Bundles and $C_0(X)$-algebras, Indiana Univ. Math. J. {\bf 45}, 463-477 (1996).

%
%
%
%
%
%

%
%

%
%

%
%

%

%
%

%
%
%

\bibitem{Ste}
N. Steenrod: {\it Topology of Fibre Bundles}, Princeton Mathematical Series,
Princeton University Press, 1965.


\bibitem{SW71}
R.F. Streater, I.F. Wilde: Fermion states of a boson field, Nucl. Phys.B {\bf 24}, 561-575 (1970).

%
%

\bibitem{Vas03}
E. Vasselli: Crossed Products by Endomorphisms, Vector Bundles and Group Duality, Int.J.Math. {\bf 16}(2), 137-171 (2005).


\bibitem{Vas05}
 E. Vasselli: Crossed Products by endomorphisms, Vector Bundles and Group Duality, II,  Int. J. Math. 17(1) (2006) 65-96.

%
%

\end{thebibliography}


}
\end{document}